\documentclass[11pt]{amsart}
\usepackage{amstext,amssymb,amsmath,amsbsy,dsfont,tikz}

 \usepackage[left=2.65cm,right=2.65cm,top=3.2cm,bottom=3.2cm]{geometry}

\usepackage{amstext,amssymb,amsmath,amsbsy}

\usepackage{tikz}
\usepackage{hyperref,cleveref}
\usepackage{amscd}
\usepackage{amsfonts}
\usepackage{indentfirst}   
\usepackage{verbatim}
\usepackage{amsmath}
\usepackage{amsthm}
\usepackage{enumerate}
\usepackage{graphicx}
\usepackage{subfig}
\usepackage{color}
\usepackage[utf8]{inputenc}
\usepackage[T1]{fontenc}     
\usepackage[english]{babel} 
\usepackage{amssymb}
\usepackage{bm}
\usepackage{ifthen, soulutf8}
\usepackage[toc,page]{appendix}
\usepackage{mathtools}
\usepackage{etoolbox}

\DeclareFontFamily{U}{matha}{\hyphenchar\font45}
\DeclareFontShape{U}{matha}{m}{n}{
<-6> matha5 <6-7> matha6 <7-8> matha7
<8-9> matha8 <9-10> matha9
<10-12> matha10 <12-> matha12
}{}
\DeclareSymbolFont{matha}{U}{matha}{m}{n}

\DeclareFontFamily{U}{mathx}{\hyphenchar\font45}
\DeclareFontShape{U}{mathx}{m}{n}{
<-6> mathx5 <6-7> mathx6 <7-8> mathx7
<8-9> mathx8 <9-10> mathx9
<10-12> mathx10 <12-> mathx12
}{}
\DeclareSymbolFont{mathx}{U}{mathx}{m}{n}

\DeclareMathDelimiter{\vvvert} {0}{matha}{"7E}{mathx}{"17}%

\DeclarePairedDelimiterX{\normiii}[1]
{\vvvert}
{\vvvert}
{\ifblank{#1}{\:\cdot\:}{#1}}

\newtheorem{theorem}{Theorem}
\newtheorem{lem}{Lemma}
\newtheorem{lemma}{Lemma}
\newtheorem{cor}{Corollary}

\newtheorem{prop}{Proposition}
\newtheorem{proposition}{Proposition}

\newtheorem{rem}{Remark}
\newtheorem{remark}{Remark}

\newtheorem{definition}{Definition}
      \newcommand{\mC}{\mathbb{C}}

\numberwithin{equation}{section}
\numberwithin{theorem}{section}
\numberwithin{remark}{section}
\numberwithin{rem}{section}
\numberwithin{prop}{section}
\numberwithin{lem}{section}
\numberwithin{definition}{section}
\numberwithin{cor}{section}
\numberwithin{lemma}{section}
  \numberwithin{proposition}{section}

\newcommand{\T}{{\rm \bold{T}}}

\newcommand{\cN}{{\mathcal N}}

\newcommand{\bT}{{\bf T}}

\newcommand{\cF}{{\mathcal F}}

\newcommand{\cT}{\mathcal T}

\newcommand{\supp}{\operatorname{supp}}
\newcommand{\dist}{\operatorname{dist}}

\newcommand{\dive}{\operatorname{div}}

\newcommand{\eps}{\varepsilon}
\renewcommand{\epsilon}{\varepsilon}

\newcommand{\loc}{_{loc}}

\newcommand{\mN}{\mathbb{N}}

\newcommand{\mR}{\mathbb{R}}
\newcommand{\R}{\mathbb{R}}

\newcommand{\C}{\mathbb{C}}

\newcommand{\hu}{\hat u}

\newcommand{\bK}{{\rm \bold{K}}}

\newcommand{\hv}{\hat v}

\newcommand{\bR}{\mathbf{R}}

\newboolean{showcomments}
\setboolean{showcomments}{true}
\ifthenelse{\boolean{showcomments}}
{ \newcommand{\mynote}[3]{
    \fbox{\bfseries\sffamily\scriptsize#1}
    {\small$\blacktriangleright$\textsf{\textit{\color{#3}{#2}}}$\blacktriangleleft$}}}
{ \newcommand{\mynote}[3]{}}

\setul{0.5ex}{0.3ex}
\setulcolor{red}

\title[The Weyl law of transmission eigenvalues]{The Weyl law of transmission eigenvalues and the completeness of generalized transmission eigenfunctions without complementing conditions}

\author{Jean Fornerod}
\address[Jean Fornerod]{Ecole Polytechnique F\'ed\'erale de Lausanne, EPFL,  SB, GR-PI,
	\newline\indent Station 8,  CH-1015 Lausanne, Switzerland.}
\email{jean.fornerod@epfl.ch}

    \author[H.-M. Nguyen]{Hoai-Minh Nguyen}
  \address[H.-M. Nguyen]{Laboratoire Jacques Louis Lions,  \newline\indent
	Sorbonne Universit\'e,  \newline\indent
	4 Place Jussieu, 75252, Paris, France}
\email{hoai-minh.nguyen@sorbonne-universite.fr}

\begin{document}
\maketitle

\begin{abstract} The transmission eigenvalue problem is a system of two second-order elliptic equations of  two unknowns  equipped with the Cauchy data on the boundary.  In this work, we establish the Weyl law for the eigenvalues  and the completeness of the generalized eigenfunctions for a system without complementing conditions, i.e.,  the two equations of  the system have the same  coefficients for the second order terms, and thus being degenerate. These coefficients are allowed to be anisotropic and are assumed  to be of class $C^2$. One of the keys of the analysis is to establish the well-posedness and the regularity in $L^p$-scale for such a system. As a result, we largely extend and rediscover  known results for which the coefficients for the second order terms are required to be isotropic and of class $C^\infty$ using a new approach. 
\end{abstract}

\medskip 



\tableofcontents

\section{Introduction}

The  transmission eigenvalue problem plays a role in the   inverse scattering theory for inhomogeneous media.  This eigenvalue  problem is connected to the injectivity of the corresponding scattering operator \cite{CM88}, \cite{Kirsch86}. Transmission eigenvalues are  related to interrogating frequencies for which there is an incident field that is  not scattered by the medium.  In the acoustic setting,  the  transmission eigenvalue problem is a system of two second-order elliptic equations  of two unknowns  equipped with the Cauchy data on the boundary.   After four decades of extensive study, the spectral properties are known to depend on a type of contrasts of the media near the boundary. Natural and interesting questions on  the interior transmission eigenvalue problem include:  the {\it discreteness} of the spectrum (see e.g. \cite{CCG10, BCH11, Sylvester12, LV12, MinhHung1,CV21}), the  {\it location} of transmission eigenvalues (see  \cite{CGH10, CL12, Vodev15, Vodev18},  and also  \cite{CMV20} for the application in time domain),  the {\it Weyl law} of transmission eigenvalues and the {\it completeness} of the generalized eigenfunctions (see e.g. \cite{LV12, LV12-A, LV15, Robbiano16}).  We refer the reader to \cite{CCH16} for a recent, and self-contained introduction  to the transmission eigenvalue problem and its applications.

 Let us  describe its  mathematical formulation. 
 Let  $\Omega$ be a bounded, simply connected, open subset of $\mR^d$ of class $C^3$ with $d\geq 2$. Let $A_1,A_2$ be two real,  symmetric matrix-valued functions,  and let  $\Sigma_1,\Sigma_2$ be two bounded  positive functions that are all defined in $\Omega$. Assume that $A_1$ and $A_2$ are uniformly elliptic,  and $\Sigma_1$ and  $\Sigma_2$ are bounded below by a positive constant in $\Omega$, i.e., for some constant $\Lambda \ge 1$,  one has,  for $\ell=1, 2$, 
\begin{equation}\label{condi1}
\Lambda^{-1} |\xi|^2\leq \langle  A_\ell (x) \xi, \xi \rangle \leq \Lambda |\xi|^2  \quad \mbox{ for all } \xi \in \mR^d, \mbox{ for a.e. }  x\in \Omega, 
\end{equation}
and
\begin{equation}\label{condi2}
\Lambda^{-1}  \leq \Sigma_\ell(x)\leq \Lambda \mbox{ for a.e. } x \in \Omega.  
\end{equation}
Here and in what follows,   $\langle \cdot, \cdot \rangle$ denotes the Euclidean scalar product in $\mC^d$ and $|\cdot|$ is the corresponding norm.  

A complex number $\lambda$ is called an eigenvalue of the transmission eigenvalue problem associated with the pairs $(A_1, \Sigma_1)$ and $(A_2, \Sigma_2)$ in $\Omega$ if there is a non-zero pair of functions $(u_1, u_2) \in [H^1(\Omega)]^2$ that satisfies the system
\begin{equation}\label{ITE}  
\left\{\begin{array}{cl}
\dive(A_1 \nabla u_1) - \lambda\Sigma_1 u_1= 0 &\text{ in}~\Omega, \\[6pt]
\dive(A_2 \nabla u_2) - \lambda\Sigma_2 u_2= 0 &\text{ in}~\Omega, \\[6pt]
u_1 =u_2, \quad  A_1 \nabla u_1\cdot \nu = A_2 \nabla u_2\cdot \nu  & \text{ on }\Gamma. 
\end{array}\right. 
\end{equation}
Here and in what follows,  $\Gamma$ denotes $\partial \Omega$, and  $\nu$ denotes the outward, normal, unit vector on $\Gamma$. Such a pair  $(u_1, u_2)$ is then called an eigenfunction pair.

Assume that $A_1, \,  A_2,  \, \Sigma_1, \, \Sigma_2$ are continuous in $\bar \Omega$, and the following conditions on the boundary $\Gamma$ hold, with $\nu = \nu(x)$ : 
\begin{equation}\label{cond1}
\langle A_2(x) \nu, \nu\rangle \langle A_2(x) \xi, \xi \rangle  - \langle A_2(x) \nu,  \xi \rangle^2 
\neq \langle A_1(x) \nu, \nu\rangle \langle A_1(x) \xi, \xi \rangle  - \langle A_1(x) \nu,  \xi \rangle^2, 
\end{equation}
for all $x \in \Gamma$ and for all $\xi \in \mR^d \setminus \{0 \}$ with  $\langle \xi, \nu \rangle =0$,  and 
\begin{equation}\label{cond2}
\big\langle  A_2(x) \nu, \nu \big\rangle \Sigma_2(x) \neq \big\langle  A_1(x) \nu, \nu \big\rangle \Sigma_1(x),   \; \forall \, x \in \Gamma. 
\end{equation}
(Q. H.) Nguyen and the second author \cite{MinhHung2} established  the Weyl law of eigenvalues and the completeness of the generalized eigenfunctions for transmission eigenvalue problem under conditions \eqref{cond1} and \eqref{cond2}  via the Fourier analysis assuming that $A_1, \,  A_2,  \, \Sigma_1, \, \Sigma_2$ are continuous in $\bar \Omega$. 
Condition \eqref{cond1} is equivalent to the celebrated complementing conditions due to Agmon, Douglis, and Nirenberg \cite{ADNII} (see also \cite{ADNI}). The explicit formula given here was derived in \cite{Ng-WP} in the context of  the study of negative index materials. Conditions  \eqref{cond1} and \eqref{cond2}  were  derived by  (Q. H.) Nguyen and the second author in \cite{MinhHung1} in their study of the discreteness of the eigenvalues for transmission eigenvalue problem.  

In the case
\begin{equation}\label{condi3}
A_1 = A_2 = A\mbox{ in $\Omega$}, 
\end{equation}
it was also shown by (Q. H.) Nguyen and the second author \cite{MinhHung1} (see also \cite{Sylvester12}) using the multiplier technique that the discreteness holds if 
\begin{equation}\label{condi4}
\Sigma_1 \neq  \Sigma_2 \mbox{ on } \Gamma. 
\end{equation}

The goal of this paper is to study the Weyl law of the eigenvalues and the completeness of the generalized eigenfunctions under conditions \eqref{condi3} and \eqref{condi4}.  It is worth noting that results in this direction have been obtained previously with more constraints on the coefficients than \eqref{condi3} and \eqref{condi4}.  Robbiano \cite{Robbiano16} (see also \cite{Robbiano13}) gives the sharp order of the counting number when $A =  I$ in $\Omega$, $\Sigma_1 = 1$, $\Sigma_2 \neq \Sigma_1$ near the boundary and $\Sigma_2$ is smooth.  The analysis is based on both the microanalysis (see,  e.g.,  \cite{GG94, Zworski12}) and the regularity theory for the transmssion eigenvalue problem. In the isotropic case, the Weyl law  was established by Petkov and Vodev \cite{PV17} and Vodev \cite{Vodev18, Vodev18-2, Vodev19} for $C^\infty$ coefficients.  Their analysis  is heavily based on  microanalysis and  the smoothness condition is strongly required. In addition, their work involved a delicate analysis on the Dirichlet to Neumann maps using non-standard parametrix construction  initiated by Vodev \cite{Vodev15}, which have their own interests.  It is not clear how one can improve the $C^\infty$ condition  and  extend their results to the anisotropic setting using their approach.  Concerning the completeness of the generalized eigenfunctions, we want to mention the work of Robbiano \cite{Robbiano13} where the case $A = I$ and $\Sigma_1 \neq  \Sigma_2$ in $\bar \Omega$ was considered. 

\medskip 
We are ready to state the main results of this paper. From now on, we will assume in addition that 
\begin{equation}\label{condi5}
\|(A_1,A_2)\|_{C^2(\bar \Omega)} + \|(\Sigma_1,\Sigma_2)\|_{C^1(\bar \Omega)} \leq \Lambda. 
\end{equation}
We  denote by $(\lambda_j)_j$ the set of transmission eigenvalues associated with the transmission eigenvalue problem \eqref{ITE}.
 
 \medskip  
   
Concerning the Weyl law, we have 
   
\begin{theorem}
\label{thm1}
Assume \eqref{condi1}-\eqref{condi2} and \eqref{condi3}-\eqref{condi5}.  Let $\cN(t)$ denote the counting function, i.e.
\[
\cN(t) = \#\{j\in \mathbb{N} : |\lambda_j|\leq t\}. 
\]
Then
\[
\mathcal{N}(t) = {\bf c}t^{\frac{d}{2}} + o(t^{\frac{d}{2}}) \mbox{ as } t \to + \infty, 
\]
where 
\[
\mathbf{c} : = \frac{1}{(2\pi)^d} \sum_{\ell =1}^2   \int_\Omega\left | \Big \{ \xi \in \mR^d;  \; \langle A_\ell(x) \xi, \xi \rangle < \Sigma_\ell (x) \Big\} \right | dx.
\]
\end{theorem}

For a measurable subset $D$ of $\mR^d$, we denote $|D|$ its  (Lebesgue) measure. 

\medskip 

Concerning the completeness, we obtain

\begin{theorem}\label{thm2}
Assume \eqref{condi1}-\eqref{condi2} and \eqref{condi3}-\eqref{condi5}. The set of generalized eigenfunction pairs of \eqref{ITE} is complete in $L^2(\Omega) \times L^2(\Omega)$.
\end{theorem}

\begin{rem} \rm As a direct consequence of either \Cref{thm1} or \Cref{thm2}, the number of eigenvalues of the transmission eigenvalue problem is infinite. As far as we know, this fact is new under the  assumption that  $A$ is allowed to be anisotropic and the regularity of the coefficients are only required up to the order 2. 
\end{rem}

\medskip 
The analysis used in the proof of \Cref{thm1} and/or \Cref{thm2} also allows us to obtain the following result on the  transmission eigenvalue free region of the complex plane ${\mathbb C}$.

\begin{proposition}\label{pro}  Assume \eqref{condi1}-\eqref{condi2} and \eqref{condi3}-\eqref{condi5}.  For $\gamma >0$, there exists $\lambda_0 > 0$ such that if $\lambda \in \mC$ with $|\Im (\lambda)| \ge \gamma  |\lambda|$ and $|\lambda| \ge \lambda_0$, then $\lambda$ is not a transmission eigenvalue. 
\end{proposition}

Here and and in what follows, for $z \in \mC$, let  $\Im(z)$ denote the imaginary part of $z$. 

\medskip 
A more general result of \Cref{pro} is given in  \Cref{pro-Free}.

\begin{rem} \rm \label{rem-cm}
Since $\gamma > 0$  can be chosen arbitrary small, combining the discreteness result in \cite{MinhHung1} mentioned above and  Proposition \ref{pro}, one derives that all the transmission eigenvalues, but finitely many, lie in a wedge of arbitrary small angle. 
\end{rem}

Some comments on \Cref{thm1} and \Cref{thm2}  are in order. In the conclusion of  \Cref{thm1}, the multiplicity of eigenvalues is taken into account and the multiplicity is associated with some operator $T_{\lambda^*}$, which is introduced in \Cref{sect-W} (see \eqref{pre-Tl} and \eqref{def-lambda^*}). Concerning $T_{\lambda^*}$, the following facts hold (see \Cref{rem-comp} and \Cref{rem-indep} for more information):  if $\mu$ is a characteristic value of the operator $T_{\lambda^*}$ associated with an eigenfunction $(u,v)$ and $\lambda^*+\mu\neq 0$, then $ \lambda^*+ \mu$ is a transmission eigenvalue of \eqref{ITE}  with an eigenfunction pair $(u_1,u_2)$ given by 
\[
u_1 = (\lambda^*+\mu)u+v \quad \quad \mbox{ and }\quad \quad u_2 = v. 
\]
Moreover, if $\lambda_j$ is a transmission eigenvalue problem, then $\lambda_j \neq \lambda^*$ and $\lambda_j - \lambda^*$ is a characteristic value of $T_{\lambda^*}$. 
In \Cref{thm2}, the generalized eigenfunctions are also associated to such an operator $T_{\lambda^*}$. We recall that the generalized eigenfunctions are complete in  $[L^2(\Omega)]^2$ if the subspace spanned by them is dense in  $[L^2(\Omega)]^2$. 

\medskip 

\Cref{thm1} and \Cref{thm2} provide the Weyl laws and the completeness under the assumptions \eqref{condi3} and \eqref{condi4} assuming the regularity conditions in \eqref{condi5}.  Our results hold for $A_1 = A_2 = A$ being anisotropic  in contrast to the isotropic setting considered previously. Moreover, the regularity assumption \eqref{condi5} on the coefficients was out of reach previously. 

\medskip 
Our approach is in the spirit of \cite{MinhHung2}
and is hence different from the ones used to study these problems given in the previous works mentioned above. The key idea is to establish the {\it regularity theory} for the transmission eigenvalue problem under the stated assumptions (see \Cref{thm-WP}). Nevertheless, several new ingredients and observations are required for the regularity theory due to the fact that \eqref{condi3}, which is degenerate,  is considered instead of \eqref{cond1}. One of the key steps to capture the phenomena is to  derive appropriate estimates in a half plane setting. It is important to note that since $A_1 = A_2 = A$, the setting is non-standard,   and the classical arguments pioneered  in \cite{ADNI, ADNII} cannot be applied since the role of $\Sigma_1$ and $\Sigma_2$ are ignored  there. To this end, our arguments for the Cauchy problems not only require the  information of the first derivatives and their structure of the data but also involve the information of the second derivatives and their structure (see, e.g., \Cref{lem-vwp}). This is quite distinct   from the complementing case where the arguments for the Cauchy problems  only require the  information of the first derivatives and no structure of the data is required \cite{MinhHung2} (see, e.g., \cite[Lemma 2 and Corollary 2]{MinhHung2}). 
One might note that the arguments used to derive the discreteness in \cite{MinhHung1} requires less assumption on the regularity of the coefficients but only give the information for one direction of $\lambda$ ($\arg{\lambda} = \pi/2$) for large $\lambda$. This is not sufficient to apply the theory of Hilbert-Schmidt operators. 

\medskip We have so far discuss the transmission eigenvalue problem in the acoustic setting. Known results for the transmission eigenvalue problem in the electromagnetic setting are much less. In this direction, we  
mention the work of Cakoni and Nguyen \cite{Cakoni-Ng21} on the state of art on the discreteness of the eigenvalues,  the work of Fornerod and Nguyen \cite{Fornerod-Ng1} on the completeness of generalized of eigenfunctions and the upper bound of the eigenvalues  for the setting considered in \cite{Cakoni-Ng21}, and the work of Vodev \cite{VodevM21} on the free region of eigenvalues for a setting considered in \cite{Cakoni-Ng21}, and the references therein. 

\medskip 
The Cauchy problem also naturally appears in the context of negative-index materials after using reflections as initiated in \cite{Ng-Complementary} (see also \cite{Ng-Superlensing-Maxwell}).  The well-posedness and the limiting absorption principle for the Helmholtz equation with sign-changing coefficients
were developed by the second author \cite{Ng-WP} using the Fourier and multiplier approach (see also \cite{NgSil}). The  work  \cite{Ng-WP} deals with the stability question of negative index materials,  and is the starting point for the analysis of the  transmission eigenvalue problems  in \cite{MinhHung1, MinhHung2} (see also \cite{Cakoni-Ng21}).  Other aspects and applications of negative-index materials as well as the stability and instability the Cauchy problem are discussed in \cite{Ng-Superlensing,  Ng-Negative-Cloaking, Ng-CALR, Ng-CALR-O} and the references therein. A survey is given in \cite{Ng-Survey}. 

\medskip 
The paper is organized as follows. \Cref{sect-Notations} is devoted to define some notations used throughout the paper.  In \Cref{sect-WP}, we establish the well-posedness and the regularity theory for the Cauchy systems associated with the transmission eigenvalue problems. The analysis is then developed in such a way that the theory of Hilbert-Schmidt operators can be used. This is given in \Cref{sect-W} where the Weyl laws are established. The completeness is considered in \Cref{sect-C}.

\section{Notations}\label{sect-Notations}

Here are some useful notations used throughput this paper.
We denote, for $\tau > 0$,  
\begin{equation}\label{def-Omegatau}
\Omega_\tau=\Big\{x\in \Omega: \dist(x,\Gamma)<\tau \Big\}.
\end{equation}
For $d \ge 2$, set
$$
\mR^d_+ = \Big\{x \in \mR^d; x_d > 0 \Big\} \quad \mbox{ and } \quad \mR^d_0 = \Big\{x \in \mR^d; x_d = 0 \Big\}.  
$$
We will identify $\mR^d_0$ with $\mR^{d-1}$ in several places. For $s>0$,  we denote
$$B_s = \{x \in \R^d : |x|<s\}.$$ For $m \ge 1$, $p \ge 1$, and $\lambda \in \C^*$ and $u \in W^{m,p}(\Omega) $, we define
\begin{equation}\label{def-norm-lambda}
\|u\|_{W^{m,p}_{\lambda}(\Omega)} = \left( \sum_{j=0}^m \| |\lambda|^{\frac{m-j}{2}} \nabla^ju\|_{L^p(\Omega)}^p \right)^{1/p}. 
\end{equation}

\section{Well-posedness and regularity theory for the transmission eigenvalue problems} \label{sect-WP}

In this section, we study the well-posedness and the regularity theory of the Cauchy problem 
\begin{equation}\label{ITE-C}  
\left\{\begin{array}{cl}
\dive(A_1 \nabla u_1) - \lambda\Sigma_1 u_1= f_1  &\text{ in}~\Omega, \\[6pt]
\dive(A_2 \nabla u_2) - \lambda\Sigma_2 u_2= f_2  &\text{ in}~\Omega, \\[6pt]
u_1 -u_2=0, \quad  (A_1 \nabla u_1- A_2 \nabla u_2)\cdot \nu=0  & \text{ on }\Gamma, 
\end{array}\right. 
\end{equation}
under the assumptions \eqref{condi1}-\eqref{condi2}, and \eqref{condi4}-\eqref{condi5}, and 
\begin{equation}\label{condi6}
A_1 = A_2 = A \mbox{ in $\Omega_\tau$}, 
\end{equation}
for some $\tau > 0$, instead of \eqref{condi3} for appropriate $\lambda \in \mC$ and $(f_1,f_2)$ in $L^p$-scale.

\medskip 
Here is the main result of this section. 

\begin{theorem} \label{thm-WP}
Assume \eqref{condi1}-\eqref{condi2}, \eqref{condi4}-\eqref{condi5}, and \eqref{condi6}. Let $1<p<+\infty$ and $\gamma \in (0, 1)$. There exist constants $\lambda_0>0$ and $C>0$ depending on $\Omega, \,  \Lambda$,  $\tau$, $p$, and  $\gamma$ such that for $\lambda \in \C$ with $|\lambda|>\lambda_0$ and $|\Im (\lambda)| \geq \gamma |\lambda|$ and for  $(f_1,f_2)\in [L^p(\Omega)]^2$, there is a unique solution $(u_1,u_2)\in [L^p(\Omega)]^2$ with $u_1 - u_2 \in W^{2,p}(\Omega)$ of the Cauchy problem \eqref{ITE-C}. Moreover, 
\begin{equation} \label{thm-WP-1}
|\lambda| \|(u_1,u_2)\|_{L^p(\Omega)} + \|u_1-u_2\|_{W^{2,p}_{\lambda}(\Omega)}  \leq C \|(f_1,f_2)\|_{L^p(\Omega)}.
\end{equation}
Assume  in addition that  $f_1-f_2 \in W^{1,p}(\Omega)$. Then $(u_1,u_2) \in [W^{1,p}(\Omega)]^2$, $u_1-u_2 \in W^{3,p}(\Omega_{\tau/2})$,  and
\begin{equation} \label{thm-WP-2}
|\lambda| \|(u_1,u_2)\|_{W^{1,p}_\lambda (\Omega)} + \|u_1-u_2\|_{W^{3,p}_{\lambda}(\Omega_{\tau/2})} \leq C \left ( |\lambda|^{1/2}\|(f_1,f_2)\|_{L^p(\Omega)}+\|f_1-f_2\|_{W^{1,p}_\lambda (\Omega)}\right ). 
\end{equation}
\end{theorem}

\begin{remark} \rm The boundary conditions must be understood as 
$$
u_1 - u_2 = 0 \mbox{ on } \Gamma \quad \mbox{ and } \quad A \nabla (u_1 - u_2) \cdot \nu = 0 \mbox{ on } \Gamma,  
$$ 
which make sense since $u_1 - u_2 \in W^{2, p}(\Omega)$. 
\end{remark}

\begin{remark} \rm In \eqref{thm-WP-2}, we only estimate $\|u_1-u_2\|_{W^{3,p}_{\lambda}(\Omega_{\tau/2})}$ not $\|u_1-u_2\|_{W^{3,p}_{\lambda}(\Omega)}$ since $f_1$ and $f_2$ are not supposed to be in $W^{1, p}(\Omega)$. Nevertheless, when $A_1 = A_2$ in $\Omega$, the estimate is also valid for $\|u_1-u_2\|_{W^{3,p}_{\lambda}(\Omega)}$. 
\end{remark}

\begin{rem} \rm As a consequence of \eqref{thm-WP-1} and the theory of regularity of elliptic equations, one derives that $(u_1,u_2) \in [W^{2,p}_{loc}(\Omega)]^2$ and for $\Omega' \Subset \Omega$ \footnote{Recall that $\Omega' \Subset \Omega$ means $\overline{ \Omega'} \subset \Omega$.},  it holds 
\[
\|(u_1,u_2)\|_{W^{2,p}_\lambda (\Omega')} \leq C \|(f_1,f_2)\|_{L^p(\Omega)}, 
\]
where $C$ depends also on $\Omega'$ (see,  e.g.,  \cite[Lemma 17.1.5]{Hor3} and  \cite[Theorem 9.11]{GilbargTrudinger}). 
\end{rem}
 
As a consequence of \Cref{thm-WP}, we obtain the following result on the free-region of the eigenvalues. 

\begin{proposition}\label{pro-Free}  Assume \eqref{condi1}-\eqref{condi2}, \eqref{condi4}-\eqref{condi5}, and \eqref{condi6}.  For $\gamma >0$, there exists $\lambda_0 > 0$ such that if $\lambda \in \mC$ with $|\Im (\lambda)| \ge \gamma  |\lambda|$ and $|\lambda| \ge \lambda_0$, then $\lambda$ is not a transmission eigenvalue. 
\end{proposition}

The rest of this section, containing two subsections, is devoted to the proof of \Cref{thm-WP}. The first one is on  the analysis in the half space. The proof of \Cref{thm-WP} is then given in the second subsection.

\subsection{Half space analysis}
Let $1< p < +\infty$. For $j=1, 2, \cdots$, and $\lambda \in \mC \setminus \{0 \}$,  we denote 
\[
\|\psi\|_{W^{j-1/p,p}_\lambda(\R^d_0)} = |\lambda|^{1/2-1/(2p)} \|\psi\|_{W^{j-1,p}_\lambda(\R^d_0)} + | \nabla^{j-1} \psi |_{W^{1-1/p,p}(\R^d_0)},
\]
where $\|\psi\|_{W^{j-1,p}_\lambda(\R^d_0)}$ is defined as in  \eqref{def-norm-lambda} with $\Omega = \mR^d_0$, and 
\[
| \psi |_{W^{1-1/p,p}(\R^d_0)}^p = \int_{\R^{d-1}} \int_{\R^{d-1}} \frac{|\psi(x')-\psi(y')|^p}{|x'-y'|^{d+p-2}} \, dx' \, dy'.
\]
By the trace theory,  there exists a positive constant $C$ depending only on $p$ and $j$ such that
\[
\| u \|_{W^{j-1/p,p}_\lambda (\R^d_0)} \leq C \|u\|_{W^{j,p}_\lambda (\R^d_+)} \mbox{ for } u \in W^{j, p} (\mR^d_+). 
\]
In fact, this inequality holds for $\lambda \in \mC$ with $|\lambda | = 1$; the general case follows by scaling. 

\medskip 
The starting point and the key ingredient of our analysis is \Cref{lem-vwp}. \Cref{lem-main} below is a special case of \Cref{lem-vwp} and is later used to derive \Cref{lem-vwp}. 

\begin{lem} \label{lem-main}
Let $A \in \R^{d \times d}$ be  a constant symmetric matrix  and let $\Sigma_1,\Sigma_2$ be two positive constants such that 
\[
\Lambda^{-1} |\xi|^2\leq \langle  A  \xi, \xi \rangle \leq \Lambda |\xi|^2 \mbox{ for all } \xi \in \mR^d, 
\]
\[
\Lambda^{-1}  \leq \Sigma_1, \Sigma_2 \leq \Lambda, \quad \mbox{ and } \quad |\Sigma_1 - \Sigma_2| \ge  \Lambda^{-1},  
\]
for some $\Lambda \ge 1$. Let $\gamma \in (0, 1)$, $1<p<+\infty$,  and  let $\varphi \in W^{2-1/p,p}(\R^d_0)$. Given $\lambda \in \mC$ with $|\lambda| \ge 1$ and $|\Im (\lambda)| \geq \gamma |\lambda|$, there exists a unique solution $(u_1,u_2) \in [L^p(\R^d_+)]^2$ with $u_1-u_2 \in W^{2,p}(\R^d_+)$ of the following Cauchy problem
\[
\left \{
\begin{array}{cl}
\dive (A\nabla u_1)-\lambda \Sigma_1 u_1 = 0 & \mbox{ in }\R^d_+, \\[6 pt] 
\dive (A\nabla u_2)-\lambda \Sigma_2 u_2 = 0 & \mbox{ in }\R^d_+, \\[6 pt] 
u_1-u_2 = \varphi , \quad A\nabla (u_1-u_2) \cdot e_d = 0 & \mbox{ on }\R_0^d.
\end{array}
\right .
\]
Moreover, 
\begin{equation}
\label{lem-main-p1}
|\lambda|\|(u_1,u_2)\|_{L^p(\R^d_+)}  + \|u_1-u_2 \|_{W^{2,p}_\lambda (\R^d_+)} \leq C \| \varphi \|_{W^{2-1/p,p}_\lambda (\R^d_0)}. 
\end{equation}
Assume  in addition that $\varphi \in W^{3-1/p,p}(\R^d_0)$. Then $(u_1,u_2) \in [W^{1,p}(\R^d_+)]^2$ with  $u_1-u_2 \in W^{3,p}(\R^d_+)$, and 
\begin{equation}
\label{lem-main-p2}
|\lambda|\|(u_1,u_2)\|_{W^{1,p}_\lambda(\R^d_+)} + \|u_1-u_2 \|_{W^{3,p}_\lambda (\R^d_+)}  \leq C \|\varphi\|_{W^{3 - 1/p, p}_\lambda (\R^d_0)}. 
\end{equation}
Here $C$ is a positive constant depending only on  $\Lambda$, $\gamma$, $p$, and $d$. 
\end{lem}

\begin{proof}
For a function $u : \R^d \to \C$ (resp. $\varphi :\R^{d-1} \to \C$) we denote by $\hat{u}$ the Fourier transform of $u$ with respect to the first $(d-1)$ variables (resp. by $\hat{\varphi}$ the Fourier transform of $\varphi$), i.e.,  for $(\xi',x_d) \in \R^{d-1}\times (0,\infty)$, 
\[
\hat{u}(\xi',x_d) = \int_{\R^{d-1}}u(x',x_d)e^{-ix'\cdot \xi'} \, dx' \quad \mbox{ and } \quad \hat{\varphi}(\xi')= \int_{\R^{d-1}}\varphi (x') e^{-ix'\cdot \xi'} \, dx' .
\]

Since, for $\ell =1,2$,
\begin{equation*}
\dive(A \nabla u_\ell) -  \lambda \Sigma_\ell u_\ell = 0 \mbox{ in } \mR^{d}_+,   
\end{equation*}
it follows that 
\begin{equation*}
a \hu_\ell''(\xi',t) + 2 i b (\xi') \hu_\ell'(\xi',t ) - (c (\xi') + \lambda \Sigma_\ell) \hu_\ell(\xi',t) = 0 \mbox{ for } t > 0,  
\end{equation*}
where
\begin{equation}\label{lem-main-abc}
a = \langle A e_d, e_d \rangle, \; \;  b(\xi') = \sum_{j=1}^{d-1} A_{jd} \xi'_j ,    \; \; c (\xi') = \sum_{i,j=1}^{d-1}A_{ij}\xi'_i \xi'_j, \; \;   \mbox{ and }  \; \;  a  c (\xi')  - b(\xi')^2  > 0,
\end{equation}
since $A$ is symmetric and positive.  One then obtains, see,  e.g., \cite[proof of Lemma 2]{MinhHung2} for the details, 
\begin{equation}\label{lem-main-ul}
\quad \hu_\ell(\xi', t) = \alpha_\ell (\xi') e^{\eta_\ell (\xi') t} 
\end{equation}
where 
\begin{equation}\label{lem-main-etal} 
\eta_\ell (\xi') = \frac{1}{a} \big( - i b (\xi') - \sqrt{\Delta_\ell (\xi')} \big) 
\end{equation}
and
\begin{equation}\label{lem-main-alphal}
\alpha_\ell  (\xi') = \frac{\hat{\varphi}(\xi') \sqrt{\Delta_{\ell+1} (\xi')}}{\sqrt{\Delta_2(\xi')} - \sqrt{\Delta_1 (\xi')}} \quad \mbox{ with } \quad \Delta_\ell  (\xi')=  -b^2 (\xi') + a \big(c (\xi') + \lambda \Sigma_\ell \big). 
\end{equation}
Here we use the convention $\Delta_{2 + \ell} = \Delta_\ell$, and  $\sqrt{\Delta_\ell}$ denotes the square root of $\Delta_\ell$ with the positive real part. 

Let $v_\ell \in W^{1, p} (\mR^d_+)$ for $\ell = 1, \, 2$  be the unique solution of the system 
\[
\left \{
\begin{array}{cl}
\dive (A\nabla v_\ell)-\lambda \Sigma_\ell v_\ell = 0 & \mbox{ in }\R^d_+, \\[6 pt] 
v_\ell = \varphi & \mbox{ on } \R^d_0.
\end{array}
\right .
\]
We have\footnote{The results hold for $|\lambda| =1$, see, e.g. \cite[Theorem 14.1]{ADNI},  the general case follows by scaling.} , for $\ell = 1, \, 2$,
\begin{equation}
\label{lem-main-14-1}
\|v_\ell \|_{W^{j, p}_\lambda(\R^d_+)} \leq C \| \varphi \|_{W^{j - 1/p, p}_\lambda (\R^d_0)} \mbox{ for } j =  2, 3,  
\end{equation}
and 
\begin{equation}
\label{lem-main-14-2}
\hv_\ell (\xi', t) =  \hat{\varphi}(\xi') e^{\eta_\ell (\xi') t}. 
\end{equation}

Extend $u_\ell(x', t)$ and  $\partial_{tt}^2 v_\ell(x', t)$ by 0 for  $t < 0$ for $\ell = 1, \, 2$ and {\it still} denote these extensions by $u_\ell(x', t)$ and  $\partial_{tt}^2 v_\ell(x', t)$.  Let $\cF$ denote the Fourier transform in $\mathbb{R}^d$.
 We then obtain from \eqref{lem-main-ul} and \eqref{lem-main-14-2} that, with $\xi=(\xi',\xi_d) \in \R^d$,  
$$
\cF u_\ell (\xi) =  - \frac{\hat \varphi (\xi')}{\eta_\ell (\xi')  - i \xi_d}  \frac{\sqrt{\Delta_{\ell+1} (\xi')}}{\sqrt{\Delta_2 (\xi')} - \sqrt{\Delta_1(\xi')}}  \quad \mbox{ and } \quad \cF \partial_{tt}^2 v_\ell (\xi) =  -   \frac{\hat \varphi (\xi') \eta_\ell^2 (\xi') }{\eta_\ell (\xi') - i \xi_d} . 
$$
It follows that 
\[
\mathcal{F}u_\ell (\xi) = m_{\ell,\lambda} (\xi)  \mathcal{F}\partial^2_{tt}v_\ell (\xi),
\]
where
\begin{equation}\label{def-ml}
m_{\ell,\lambda}(\xi) =  \frac{\sqrt{\Delta_{\ell+1}(\xi')}}{\eta_\ell^2 (\xi')  (\sqrt{\Delta_2(\xi')}-\sqrt{\Delta_1(\xi')})}   .
\end{equation}

Note that
\[
\Delta_2 (\xi') -\Delta_1 (\xi') = a \lambda (\Sigma_2 - \Sigma_1) \neq 0
\]
and 
\begin{multline*}
\frac{1}{\eta_\ell(\xi')} \mathop{=}^{\eqref{lem-main-etal}} \frac{a}{-i b(\xi') - \sqrt{\Delta_\ell(\xi')}} = \frac{a \big(-i b(\xi') + \sqrt{\Delta_\ell(\xi')} \big)}{-b(\xi')^2 - \Delta_\ell(\xi')}\\[6pt] 
\mathop{=}^{\eqref{lem-main-alphal}}  \frac{a \big(-i b(\xi') + \sqrt{\Delta_\ell(\xi')} \big)}{- a \big( c(\xi') + \lambda \Sigma_\ell \big)}  =  \frac{i b(\xi') - \sqrt{\Delta_\ell(\xi')}}{ c(\xi') + \lambda \Sigma_\ell}. 
\end{multline*}
We derive from  \eqref{def-ml} that 
\begin{equation}\label{def-ml1}
m_{\ell,\lambda}(\xi) = \frac{\sqrt{\Delta_{\ell+1}(\xi')}(\sqrt{\Delta_1(\xi')} + \sqrt{\Delta_2(\xi')})(ib(\xi') - \sqrt{\Delta_\ell(\xi')})^2}{a\lambda (\Sigma_2-\Sigma_1)(c(\xi') + \lambda \Sigma_\ell)^2 }. 
\end{equation}
We have, by \eqref{lem-main-abc} and  \eqref{lem-main-alphal},\footnote{Given two functions, $p_1(\xi',\lambda)$ and $p_{2}(\xi',\lambda)$ the notation $p_1(\xi,\lambda) \sim p_2(\xi',\lambda)$ means that there exists a constant $C\geq 1$ independent of $\xi'$ and $\lambda$ such that $C^{-1}|p_1(\xi',\lambda)|\leq |p_2(\xi',\lambda)|\leq C |p_1(\xi',\lambda)|$. }
$$
|\Delta_{\ell}(\xi')| \sim (|\xi'|^2 + |\lambda|),  \quad |b(\xi')| \le C |\xi'|, \quad \mbox{ and } \quad |c(\xi') + \lambda \Sigma_\ell| \sim |\xi'|^2 + |\lambda|. 
$$
We then derive from \eqref{def-ml1} that 
\begin{equation}
|\xi|^j |\nabla^j m_{\ell,\lambda}(\xi)| \le C_j |\lambda|^{-1}  \mbox{ for } j \in \mathbb{N}. 
\end{equation}
It follows from Mikhlin-H\"ormander's multiplier theorem, see,  e.g.,  \cite[Theorem 7.9.5]{Hor1}, that 
\begin{equation}\label{lem-main-P1}
|\lambda| \| u_\ell \|_{L^p(\mR^d)} \le C  \| \partial_{tt}^2 v_\ell \|_{L^p(\mR^d)},  
\end{equation}
which implies
\begin{equation}\label{lem-main-p1-1}
|\lambda| \| u_\ell \|_{L^p(\mR^d)} \mathop{\le}^{\eqref{lem-main-14-1}} C \| \varphi \|_{W^{2-1/p, p}_\lambda(\mR^d)}.   
\end{equation}

On the other hand, one has 
\begin{equation*}\left\{
\begin{array}{c}
\dive \big(A \nabla (u_1 - u_2) \big) - \lambda \Sigma_1 (u_1 - u_2) = \lambda (\Sigma_1 - \Sigma_2) u_2 \mbox{ in } \mR^d_+, \\[6pt]
u_1 - u_2 = 0 \mbox{ on } \mR^d_0. 
\end{array}\right. 
\end{equation*}
This yields 
\[
\| u_1 - u_2 \|_{W^{2, p}_\lambda(\mR^d_+)} \le C \| \lambda (\Sigma_1 - \Sigma_2) u_2 \|_{L^p(\mR^d_+)} \mathop{\le}^{\eqref{lem-main-p1-1}} C \| \varphi \|_{W^{2 - 1/p, p }_\lambda (\mR^d_+)}. 
\]

We next deal with \eqref{lem-main-p2}. By taking the derivative of the system with respect to $x_j$ for $1 \le j \le d -1$ and applying 
\eqref{lem-main-p1}, we have, for $1 \le j \le d - 1$,  
\begin{equation}
\label{lem-main-p2-1}
|\lambda| \| (\partial_{x_j} u_1, \partial_{x_j}u_2)\|_{L^p(\mR^d_+)} + \|(\partial_{x_j}u_1- \partial_{x_j} u_2)\|_{W^{2,p}_{\lambda}(\mR^d_+)}  \leq C \| \partial_{x_j} \varphi\|_{W^{2- 1/p, p}_\lambda(\mR^d_0)}.
\end{equation}
Extend $\partial_t u_\ell(x', t)$ and  $\partial_{ttt}^3 v_\ell(x', t)$ by 0 for  $t < 0$ for $\ell = 1, \, 2$ and {\it still} denote these extensions by $\partial_t u_\ell(x', t)$ and  $\partial_{ttt}^3 v_\ell(x', t)$.  
 We then obtain from \eqref{lem-main-ul} and \eqref{lem-main-14-2} that, with $\xi=(\xi',\xi_d) \in \R^d$,  
$$
\cF \partial_t u_\ell (\xi) =  - \frac{\hat \varphi (\xi') \eta_\ell(\xi')}{\eta_\ell (\xi')  - i \xi_d}  \frac{\sqrt{\Delta_{\ell+1} (\xi')}}{\sqrt{\Delta_2 (\xi')} - \sqrt{\Delta_1(\xi')}}  \quad \mbox{ and } \quad \cF \partial_{ttt}^3 v_\ell (\xi) =  -   \frac{\hat \varphi (\xi') \eta_\ell^3 (\xi') }{\eta_\ell (\xi') - i \xi_d} . 
$$
This yields 
\[
\mathcal{F} \partial_t u_\ell (\xi) = m_{\ell,\lambda} (\xi)  \mathcal{F}\partial^3_{ttt}v_\ell (\xi). 
\]
As in the proof of  \eqref{lem-main-P1}, we obtain 
\begin{equation*}
|\lambda| \| \partial_t u_\ell \|_{L^p(\mR^d)} \le C  \| \partial_{ttt}^3 v_\ell \|_{L^p(\mR^d)},  
\end{equation*}
which implies
\begin{equation}\label{lem-main-p2-2}
|\lambda| \| \partial_t u_\ell \|_{L^p(\mR^d)} \mathop{\le}^{\eqref{lem-main-14-1}} C \| \varphi \|_{W^{3-1/p, p}_\lambda(\mR^d)}.   
\end{equation}

Combining \eqref{lem-main-p2-1} and \eqref{lem-main-p2-2}, we derive that  
\begin{equation}
\label{lem-main-p2-3}
|\lambda| \| (\nabla  u_1, \nabla u_2)\|_{L^p(\mR^d_+)} + \|\nabla (u_1-  u_2)\|_{W^{2,p}_{\lambda}(\mR^d_+)}  \leq C \| \varphi\|_{W^{3- 1/p, p}_\lambda(\mR^d_0)}.
\end{equation}
Assertion \eqref{lem-main-p2} now follows from \eqref{lem-main-p2-3} and \eqref{lem-main-p1}.  The proof is complete. 
\end{proof}

We now state and prove a more general version of \Cref{lem-main}, which is the main ingredient of the proof of \Cref{thm-WP}. 

\begin{lem} \label{lem-vwp}
Let $A \in \R^{d \times d}$ be  a constant symmetric matrix  and let $\Sigma_1,\Sigma_2$ be two positive constants such that 
\[
\Lambda^{-1} |\xi|^2\leq \langle  A  \xi, \xi \rangle \leq \Lambda |\xi|^2 \mbox{ for all } \xi \in \mR^d, 
\]
and 
\[
\Lambda^{-1}  \leq \Sigma_1, \Sigma_2 \leq \Lambda, \quad \mbox{ and } \quad |\Sigma_1 - \Sigma_2| \ge  \Lambda^{-1},  
\]
for some $\Lambda \ge 1$. Let $\gamma \in (0, 1)$, $1<p<+\infty$,  and  let $f_1, f_2 \in L^p(\R^d_+)$, $G_1, G_2\in [L^p(\R^d_+)]^d$ with $G_1-G_2 \in [W^{1,p}(\R^d_+)]^d$, $\varphi \in W^{2-1/p,p}(\R^d_0)$, $\psi \in W^{1-1/p,p}(\R^d_0)$,  and let $r_{1}^{(ij)}, r_{2}^{(ij)} \in L^p(\R^d_+)$  with $r_{1}^{(ij)} - r_{2}^{(ij)} \in W^{2,p}(\R^d_+)$ for $1 \le i, j \le d$. Given $\lambda \in \mC$ with $|\lambda| \ge 1$ and $|\Im (\lambda)| \geq \gamma |\lambda|$, there exists a unique solution $(u_1,u_2) \in [L^p(\R^d_+)]^2$ with $u_1-u_2 \in W^{2,p}(\R^d_+)$ of the following Cauchy problem
\begin{equation} \label{lem-vwp-3}
\left \{
\begin{array}{cl}
\dive (A\nabla u_1)-\lambda \Sigma_1 u_1 = f_1 + \dive (G_1) + \sum_{i,j=1}^d\partial_{ij}^2 r_{1}^{(ij)}& \mbox{ in }\R^d_+, \\[6 pt] 
\dive (A\nabla u_2)-\lambda \Sigma_2 u_2 = f_2 + \dive (G_2) + \sum_{i,j=1}^d \partial_{ij}^2 r_{2}^{(ij)} & \mbox{ in }\R^d_+, \\[6 pt] 
u_1-u_2 =\varphi, \quad A\nabla (u_1-u_2) \cdot e_d = \psi & \mbox{ on }\R_0^d.
\end{array}
\right .
\end{equation}
Moreover, 
\begin{align}
\label{lem-vwp-4}
C &  \left(|\lambda|\|(u_1,u_2)\|_{L^p(\R^d_+)}  + \|u_1-u_2 \|_{W^{2,p}_\lambda (\R^d_+)} \right) \\[6pt]
 & \leq  \|(f_1,f_2)\|_{L^p(\R^d_+)}  + |\lambda|^{1/2}\|(G_1,G_2)\|_{L^p(\R^d_+)}  +  \sum_{i,j=1}^d |\lambda|\|(r_{1}^{(ij)}, r_{2}^{(ij)})\|_{L^p(\R^d_+)} \nonumber  \\[6pt]
& + \|\varphi\|_{W^{2-1/p,p}_\lambda (\R^d_0)}+\|\psi\|_{W^{1-1/p,p}_{\lambda}(\R^d_0)} + \|G_1-G_2 \|_{W^{1,p}_\lambda (\R_+^d)} 
 + \sum_{i,j=1}^d  \|r_{1}^{(ij)} -r_{2}^{(ij)}\|_{W^{2,p}_{\lambda}(\R^d_+)} \nonumber. 
\end{align}
Assume in addition  that  $f_1-f_2 \in W^{1,p}(\R^d_+)$, $G_1-G_2 \in W^{2,p}(\R^d_+)$, $\varphi \in W^{3-1/p,p}(\R^d_0)$,  $\psi \in W^{2-1/p,p}(\R^d_0)$, and $r_{1}^{(ij)} = r_{2}^{(ij)} = 0$ for all $1\leq i,j \leq d$.  Then $(u_1,u_2) \in W^{1,p}(\R^d_+)$ with $u_1-u_2 \in W^{3,p}(\R^d_+)$,  and it holds 
\begin{multline}
\label{lem-vwp-4bis}
C \left( |\lambda|\|(u_1,u_2)\|_{W^{1,p}_\lambda (\R^d_+)} + \|u_1-u_2 \|_{W^{3,p}_\lambda (\R^d_+)} \right) \\[6pt]
 \leq  |\lambda|^{1/2}\|(f_1,f_2)\|_{L^p(\R^d_+)}  + |\lambda|\|(G_1,G_2)\|_{L^p(\R^d_+)}  + \|\varphi\|_{W^{3-1/p,p}_\lambda (\R^d_0)} \\[6pt] 
+\|\psi\|_{W^{2-1/p,p}_{\lambda}(\R^d_0)} +  \|f_1-f_2\|_{W^{1,p}_\lambda (\R^d_+)}  + \|G_1-G_2 \|_{W^{2,p}_\lambda (\R_+^d)} . 
\end{multline}
Here $C$ denotes a positive constant  depending only on  $\Lambda$, $\gamma$, $d$, and $p$.  
\end{lem}

\begin{remark} \rm Concerning \eqref{lem-vwp-4bis}, the assumption $r_{1}^{(ij)} = r_{2}^{(ij)} = 0$ for all $1\leq i,j \leq d$ is just to avoid the redundancy; the same estimate holds for the appropriate assumptions on $r_\ell^{(ij)}$ but this can be put into the conditions of $f_\ell$ and $G_\ell$ instead.  
\end{remark}

\begin{proof}  Since the problem is linear, \eqref{lem-vwp-4} and \eqref{lem-vwp-4bis}  follow from the corresponding estimates in the following two cases: 

\begin{enumerate}
\item[] $\bullet$ Case 1:  $f_1=f_2 = 0$, $G_1 = G_2 = 0$,  and $r_{1}^{(ij)} = r_{2}^{(ij)} = 0$ for all $1\leq i, j \leq d$.  

\item[] $\bullet$ Case 2: $\varphi =0$ and $\psi=0$.
\end{enumerate}

 We now proceed the proof for these cases.

\medskip  
\noindent \underline{\it Case 1}:  $f_1=f_2 = 0$, $G_1 = G_2 = 0$, and $r_{1}^{(ij)}  = r_{2}^{(ij)}= 0$ for all $1\leq i, j \leq d$. We have 
\[
\left \{
\begin{array}{cl}
\dive (A\nabla u_1)-\lambda \Sigma_1 u_1 = 0& \mbox{ in }\R^d_+, \\[6 pt] 
\dive (A\nabla u_2)-\lambda \Sigma_2 u_2 = 0& \mbox{ in }\R^d_+, \\[6 pt] 
u_1-u_2 =\varphi, \quad A\nabla (u_1-u_2) \cdot e_d = \psi & \mbox{ on }\R_0^d.
\end{array}
\right .
\]

Let $v \in W^{1,p}(\R^d_+)$ be the unique solution of 
\[
\left \{
\begin{array}{rc}
\dive (A\nabla v)-\lambda \Sigma_1 v = 0 & \mbox{ in }\R^d_+, \\[6 pt] 
A\nabla v \cdot e_d =  \psi & \mbox{ on }\R^d_0. 
\end{array}
\right .
\]
As a consequence of \cite[Theorem 2.3.2.7]{Grisvard} and a scaling argument, we have 
\begin{equation}
\label{est-5}
\|v\|_{W^{2,p}_\lambda(\R^d_+)} \leq C\|\psi\|_{W^{1-1/p,p}_\lambda(\R^d_0)} \quad \mbox{ and } \quad  \|v\|_{W^{3,p}_\lambda(\R^d_+)} \leq C\|\psi\|_{W^{2-1/p,p}_\lambda(\R^d_0)}. 
\end{equation}
By the trace theory, it follows that 
\begin{equation}
\label{est-5-bis}
\|v\|_{W^{2-1/p,p}_\lambda(\R^d_0)} \leq C\|\psi\|_{W^{1-1/p,p}_\lambda(\R^d_0)} \quad \mbox{ and } \quad  \|v\|_{W^{3-1/p,p}_\lambda(\R^d_0)} \leq C\|\psi\|_{W^{2-1/p,p}_\lambda(\R^d_0)}. 
\end{equation}

Considering the system of $(u_1-v,u_2)$ and using \eqref{est-5}, and \eqref{est-5-bis}, the conclusion of this case follows from   \Cref{lem-main}.

\medskip 
\noindent \underline{\it Case 2:} $\varphi =0$, $\psi=0$. In this case, we have
\[
\left \{
\begin{array}{cl}
\dive (A\nabla u_1)-\lambda \Sigma_1 u_1 = f_1 + \dive (G_1)+ \sum_{i,j=1}^d \partial_{ij}^2 r_{1}^{(ij)} & \mbox{ in }\R^d_+, \\[6 pt] 
\dive (A\nabla u_2)-\lambda \Sigma_2 u_2 = f_2 + \dive (G_2)+ \sum_{i,j=1}^d \partial_{ij}^2 r_{2}^{(ij)} & \mbox{ in }\R^d_+, \\[6 pt] 
u_1-u_2 =0, \quad A\nabla (u_1-u_2) \cdot e_d = 0 & \mbox{ on }\R_0^d.
\end{array}
\right .
\]
For $\ell =1,2$, consider the following systems
\[
\left \{
\begin{array}{rll}
\dive (A\nabla v_{\ell}^{(0)})-\lambda \Sigma_\ell v_{\ell}^{(0)}  &= f_\ell   &\mbox{ in }\R^d_+,  \\[6 pt] 
A \nabla v_{\ell}^{(0)} \cdot e_d &= 0  &\mbox{ on }\R^d_0,
\end{array}
\right .
\quad \quad \phantom{(1\leq j \leq d)}
\]
\[
\left \{
\begin{array}{rll}
\dive (A\nabla v_{\ell}^{(j)})-\lambda \Sigma_\ell v_{\ell}^{(j)} &= (G_\ell)_j    &\mbox{ in }\R^d_+,  \\[6 pt] 
A \nabla v_{\ell}^{(j)} \cdot e_d &= 0  &\mbox{ on }\R^d_0
\end{array}
\right .
\quad \quad (1\leq j \leq d),
\]
where $(G_\ell)_j$ denotes the $j$-th component of $G_\ell$, and
\[
\left \{
\begin{array}{rll}
\dive (A\nabla v_{\ell}^{(ij)})-\lambda \Sigma_\ell v_{\ell}^{(ij)} &= r_{\ell}^{(ij)}  &\mbox{ in }\R^d_+,  \\[6 pt] 
A \nabla v_{\ell}^{(ij)} \cdot e_d &= 0  &\mbox{ on }\R^d_0
\end{array}
\right .
\quad \quad (1\leq i,j \leq d).
\]
We have, see,  e.g.,  \cite[Theorem 14.1]{ADNI}, for $1\leq i,j\leq d$,
\begin{equation}
\label{lem-vw-2}
\left\{\begin{array}{c}
\|v_{\ell}^{(0)} \|_{W^{2,p}_\lambda (\R^d_+)} \leq C \|f_\ell\|_{L^p(\R^d_+)},  \\[6pt]
\|v_{\ell}^{(j)} \|_{W^{2,p}_\lambda (\R^d_+)}   \leq C \|G_\ell \|_{L^p(\R^d_+)}, \\[6pt]  
\|v_{\ell}^{(ij)} \|_{W^{2,p}_\lambda (\R^d_+)}  \leq C \| r_{\ell}^{(ij)} \|_{L^p(\R^d_+)}. 
\end{array}\right.
\end{equation}
Since, we have
\[
\left \{
\begin{array}{rll}
\dive (A\nabla( v_{1}^{(0)}-v_{2}^{(0)}))-\lambda \Sigma_1 (v_{1}^{(0)}-v_{2}^{(0)})  &= f_1-f_2 + \lambda (\Sigma_1 -\Sigma_2) v_2^{(0)}   &\mbox{ in }\R^d_+,  \\[6 pt] 
A \nabla (v_{1}^{(0)}-v_{2}^{(0)}) \cdot e_d &= 0  &\mbox{ on }\R^d_0,
\end{array}
\right .
\]
and the equations for $v_1^{(j)} - v_2^{(j)}$ and $v_1^{(ij)} - v_2^{(ij)} $ are similar, 
we also get, for $1\leq i,j\leq d$,  by using \eqref{lem-vw-2},
\begin{equation}\label{lem-vw-cc1}
\left\{\begin{array}{c}
C\| v_1^{(0)} - v_2^{(0)} \|_{W^{2, p}_\lambda(\mR^d_+)} \le  \| (f_1, f_2 )\|_{L^p(\mR^d_+)}, \\[6pt] 
C\| v_1^{(j)} - v_2^{(j)} \|_{W^{3, p}_\lambda(\mR^d_+)} \le  \| G_1 - G_2 \|_{W^{1, p}_\lambda(\mR^d_+)}+|\lambda|^{1/2}\|G_2\|_{L^p(\R^d_+)}, \\[6pt] 
C\| v_1^{(ij)} - v_2^{(ij)} \|_{W^{4, p}_\lambda(\mR^d_+)} \le  \| r_1^{(ij)} -  r_2^{(ij)}  \|_{W^{2, p}_\lambda(\mR^d_+)}+|\lambda|\|r_2^{(ij)}\|_{L^p(\R^d_+)}, 
\end{array}\right.
\end{equation}
and 
\begin{equation}\label{lem-vw-cc1-bis}
\begin{aligned}
C\| v_1^{(0)} - v_2^{(0)} \|_{W^{3, p}_\lambda(\mR^d_+)} &\le  \|f_1 - f_2 \|_{W^{1, p}_\lambda(\mR^d_+)}+|\lambda|^{1/2}\|f_2\|_{L^p(\R^d_+)}, \\
 C\| v_1^{(j)} - v_2^{(j)} \|_{W^{4, p}_\lambda(\mR^d_+)} &\le  \| G_1 - G_2 \|_{W^{2, p}_\lambda(\mR^d_+)}+|\lambda|\|G_2\|_{L^p(\R^d_+)}. 
 \end{aligned}
\end{equation}

For $\ell=1,2$, set
\[
w_\ell = v_{\ell}^{(0)} + \sum_{j=1}^d \partial_j v_{\ell}^{(j)} + \sum_{i,j=1}^d \partial^2_{ij} v_{\ell}^{(ij)}.
\]
We have
\[
\dive (A\nabla w_\ell) - \lambda \Sigma_\ell w_\ell = f_\ell+ \dive (G_\ell) + \sum_{i,j=1}^d \partial_{ij}^2 r_{\ell}^{(ij)} \mbox{ in }\R^d_+.
\] 
Moreover, 
\begin{multline}
\label{est-0}
C|\lambda|\|(w_1,w_2)\|_{L^p(\R^d_+)} \\[6pt]
\mathop{\leq}^{\eqref{lem-vw-2}} \|(f_1,f_2)\|_{L^p(\R^d_+)} + |\lambda|^{1/2}\|(G_1,G_2)\|_{L^p(\R^d_+)}+|\lambda|\sum_{i,j=1}^d\|(r_{1}^{(ij)}, r_{2}^{(ij)})\|_{L^p(\R^d_+)}.
\end{multline}
Using  \eqref{lem-vw-cc1}  and the trace theory, we  derive that
\begin{multline}
\label{est-3}
 \|w_1 - w_2 \|_{W^{2, p}_\lambda(\mR^d_+)}  +  \|w_1-w_2\|_{W^{2-1/p,p}_{\lambda}(\R^d_0)} + \|A\nabla (w_1-w_2) \cdot e_d\|_{W^{1-1/p,p}_{\lambda}(\R^d_0)} \\[6pt]
\leq C \Big( \|f_2\|_{L^p(\R^d_+)}+|\lambda|^{1/2} \|(G_1,G_2)\|_{L^p(\R^d_+)} + \sum_{i,j=1}^d |\lambda| \|r_2^{ij}\|_{L^p(\R^d_+)} \\
+\|f_1 -  f_2\|_{L^p(\R^d_+)}  
+ \|G_1-G_2\|_{W^{1,p}_{\lambda}(\R^d_+)}+\sum_{i,j=1}^d \|r_{1}^{(ij)} - r_{2}^{(ij)}\|_{W^{2,p}_\lambda (\R^d_+)} \Big ).  
\end{multline}
Considering the system of $(u_1-w_1,u_2-w_2)$,  and using \eqref{est-0} and \eqref{est-3}, assertion \eqref{lem-vwp-4} now follows from  case 1. 

To deal with assertion \eqref{lem-vwp-4bis}, instead of \eqref{est-0} and \eqref{est-3},  we use, since $r_1^{(ij)} = r_2^{(ij)} = 0$, 
\begin{equation}
\label{est-0-bis}
|\lambda|\|(w_1,w_2)\|_{W^{1, p}_\lambda(\R^d_+)} 
\mathop{\leq}^{\eqref{lem-vw-2}} C \Big( |\lambda|^{1/2} \|(f_1,f_2)\|_{L^p(\R^d_+)} + |\lambda| \|(G_1,G_2)\|_{L^p(\R^d_+)} \Big), 
\end{equation}
and 
\begin{multline}
\label{est-3-bis}
\|w_1 - w_2 \|_{W^{3, p}_\lambda(\mR^d_+)} +  \|w_1-w_2\|_{W^{3-1/p,p}_{\lambda}(\R^d_0)} + \|A\nabla (w_1-w_2) \cdot e_d\|_{W^{2-1/p,p}_{\lambda}(\R^d_0)}   \\[6pt]
\mathop{\leq}^{\eqref{lem-vw-cc1-bis}} C \Big( \|(f_1 -  f_2)\|_{W^{1, p}_\lambda(\R^d_+)} 
+ \|G_1-G_2\|_{W^{2,p}_{\lambda}(\R^d_+)} \Big). 
\end{multline}
By  considering the system of $(u_1-w_1,u_2-w_2)$, assertion \eqref{lem-vwp-4bis} now follows from  case 1.

\medskip 
The proof is complete. 
\end{proof}

\subsection{Proof of \Cref{thm-WP}} The proof is divided into two steps: 

\begin{itemize}

\item Step 1: Assuming the solution exists, we establish \eqref{thm-WP-1} and \eqref{thm-WP-2}. 

\item Step 2: We establish the existence of the solutions. 

\end{itemize}

We now proceed these two steps. 

\medskip 
\noindent {\underline{\it Step 1:}}
 For $(f_1,f_2) \in [L^p(\Omega)]^2$, let $(u_1,u_2) \in [L^p(\Omega)]^2$ with $u_1-u_2 \in W^{2,p}(\Omega)$ be a solution of \eqref{ITE-C}. We prove that \eqref{thm-WP-1} and \eqref{thm-WP-2} hold.

Applying \Cref{lem-vwp} and the freezing coefficient technique,  we deduce that there exists $\tau_* \in (0, \tau/2)$ depending only on $\Omega$, $\Lambda$, $\tau$,  and $p$, such that   
\begin{multline}
\label{thm-WP-pr1}
C  \left ( |\lambda| \|(u_1,u_2) \|_{L^p(\Omega_{\tau_*})}+\|u_1-u_2\|_{W^{2,p}_\lambda (\Omega_{\tau_*})} \right ) \\[6 pt]
\leq  \|(f_1,f_2)\|_{L^p(\Omega)} + |\lambda|^{1/2} \| (u_1,u_2) \|_{L^p(\Omega_\tau)}+ \|u_1-u_2\|_{W^{1,p}_\lambda (\Omega_\tau)}
\end{multline} 
and 
\begin{multline}
\label{thm-WP-pr1-bis}
C  \left ( |\lambda| \|(u_1,u_2) \|_{W^{1, p}_\lambda(\Omega_{\tau_*})}+\|u_1-u_2\|_{W^{3,p}_\lambda (\Omega_{\tau_*})} \right )\\[6pt]
\leq  |\lambda|^{1/2} \|(f_1,f_2)\|_{L^p(\Omega)}  +  \|f_1 - f_2 \|_{W^{1, p}_\lambda(\Omega)} \\[6pt]
+ |\lambda| \|(u_1,u_2) \|_{L^p (\Omega_\tau)}
+ |\lambda| \|(u_1,u_2) \|_{L^p (\Omega_\tau)}
+\|u_1-u_2\|_{W^{2,p}_\lambda (\Omega_\tau)},  
\end{multline} 
for every $\lambda \in \C$ with $|\Im (\lambda)| \geq c |\lambda|$ and $|\lambda| \ge 1 $. Here and in what follows, $C$ denotes a positive constant depending only on $\Omega$, $\Lambda$, $\tau$, and $p$.

Let us emphasize here that  the terms $(r_{1,ij},r_{2,ij})$  in \Cref{lem-vwp},  play a crucial role in the proof of \eqref{thm-WP-pr1} since the solutions $(u_1, u_2)$ considered are only in $[L^p(\Omega)]^2$, but not in $[W^{1, p}(\Omega)]^2$. Indeed, let consider a small neighborhood of $x_0 \in \Gamma$. Using a change of variables, without loss of generality, one might assume that the boundary in this neighbourhood  is {\it flat} already and $A_1 = A_2 = A$ there. In the freezing process, one has, in such a neighborhood, 
\begin{multline*}
\dive (A(x_0) \nabla u_\ell) - \lambda \Sigma_\ell (x_0) u_\ell = \dive \big( (A(x_0) - A(x)) \nabla u_\ell \big) + 
\dive \big( A(x) \nabla u_\ell \big) -  \lambda \Sigma_\ell (x_0) u_\ell \\[6pt]
 =\sum_{i,j=1}^d \partial^2_{ij} \Big ( ( A_{ij}(x_0) - A_{ij}(x) )  u_\ell \Big ) -  \sum_{i,j=1}^d  \partial_i \Big ( u_\ell  \partial_{j}  (A_{ij}(x_0) - A_{ij}(x) ) \Big )   + f_\ell  +   \lambda \big( \Sigma_\ell (x) - \Sigma_\ell (x_0) \big) u_\ell. 
\end{multline*}

Let $\chi \in C^\infty (\R^d)$ with the support in a sufficiently small neighborhood of $x_0$, then with $v_\ell = \chi u_\ell$ for 
 $\ell =1,2$,  we have 
\begin{multline}\label{thm-WP-cc}
\dive (A(x_0) \nabla v_\ell) - \lambda \Sigma_\ell (x_0) v_\ell 
\\[6pt]
=  \chi \sum_{i,j=1}^d \partial^2_{ij} \Big ( ( A_{ij}(x_0) - A_{ij}(x) )  u_\ell \Big )   -  \sum_{i,j=1}^d  \chi  \partial_i \Big ( u_\ell  \partial_{j}  (A_{ij}(x_0) - A_{ij}(x) ) \Big) \\[6pt]
+ \chi  f_\ell  +    \lambda \big( \Sigma_\ell (x) - \Sigma_\ell (x_0) \big) v_\ell   -  u_\ell \dive (A (x_0) \nabla \chi ) + 2 \dive (u_\ell A (x_0) \nabla \chi). 
\end{multline}
The terms $r_{\ell, ij}$ are then $( A_{ij}(x_0) - A_{ij}(x) ) \chi u_\ell = ( A_{ij}(x_0) - A_{ij}(x) ) v_\ell $. 
Since $A_1 = A_2 = A$ in $\Omega_\tau$, $u_1 - u_2 = 0$ in $\Gamma$, 
and $A \nabla (u_1 -u_2) \cdot \nu = 0$ on $\Gamma$, it follows that 
$$
v_1 - v_2 = 0 \mbox{ on } \Gamma \quad \mbox{ and } \quad A(x_0) \nabla (v_1 - v_2)  \cdot \nu    = \chi (A(x_0)-A(x))\nabla (u_1-u_2) \cdot \nu \mbox{ on } \Gamma.  
$$ 
We are thus in the situation to apply  \Cref{lem-vwp} and the freezing coefficient technique to derive \eqref{thm-WP-pr1}. 

Concerning \eqref{thm-WP-pr1-bis}, in \eqref{thm-WP-cc}, one writes $\partial^2_{ij} \Big ( ( A_{ij}(x_0) - A_{ij}(x) )  u_\ell(x) \Big ) $ under the form 
$$
\partial_{i} \Big ( ( A_{ij}(x_0) - A_{ij}(x) )  \partial_j u_\ell \Big )  + \partial_i \Big( \partial_j ( A_{ij}(x_0) - A_{ij}(x) ) u_\ell \Big).  
$$
We are thus in the situation to apply  \Cref{lem-vwp} and the freezing coefficient technique to derive \eqref{thm-WP-pr1-bis}. 
The details of the rest of the proof of \eqref{thm-WP-pr1} and  \eqref{thm-WP-pr1-bis} are omitted.

On the other hand, since
\[
\dive (A_\ell \nabla u_\ell) -\lambda \Sigma_\ell u_\ell = f_\ell  \quad \mbox{ in }\Omega,
\] 
we have, for $|\lambda| \ge 1$,  
\begin{equation}
\label{thm-WP-pr2}
\|u_\ell\|_{W^{1, p}_\lambda(\Omega \setminus \Omega_{\tau_*/4})} \leq C  \Big( |\lambda|^{-1/2}\|f_\ell\|_{L^p(\Omega)} + \|u_\ell \|_{L^p(\Omega_{\tau_*})} \Big),
\end{equation}
and
\begin{equation}
\label{thm-WP-pr3}
\|u_\ell\|_{W^{2, p}_\lambda(\Omega \setminus \Omega_{\tau_*/2})} \leq C  \Big( \|f_\ell\|_{L^p(\Omega)} + \|u_\ell \|_{W^{1, p}_\lambda(\Omega_{\tau_*} \setminus \Omega_{\tau_*/4})} \Big).  
\end{equation}

Combining \eqref{thm-WP-pr2} and \eqref{thm-WP-pr3} yields 
\begin{equation}
\label{thm-WP-pr23}
\|u_\ell\|_{W^{2, p}_\lambda(\Omega \setminus \Omega_{\tau_*/2})} \leq C \left (\|f_\ell\|_{L^p(\Omega)} + \|u_\ell \|_{L^p(\Omega_{\tau_*})} \right ). 
\end{equation}
From \eqref{thm-WP-pr1} and  \eqref{thm-WP-pr23}, we obtain 
\begin{equation*}
|\lambda| \|(u_1,u_2) \|_{L^p(\Omega)}+\|u_1-u_2\|_{W^{2,p}_\lambda (\Omega)} \le  C  \|(f_1,f_2)\|_{L^p(\Omega)} 
\end{equation*} 
for $|\lambda| \ge \lambda_0$ and for $\lambda_0$ large enough. This completes the proof of  \eqref{thm-WP-1}. 

From \eqref{thm-WP-pr1-bis}, \eqref{thm-WP-pr2}, after using  \eqref{thm-WP-1}, we obtain 
\begin{equation*}
|\lambda| \|(u_1,u_2) \|_{W^{1, p}_\lambda(\Omega)}+\|u_1-u_2\|_{W^{3,p}_\lambda (\Omega_{\tau_*})} \le  C |\lambda|^{1/2} \|(f_1,f_2)\|_{L^p(\Omega)} +  \|f_1 - f_2\|_{W^{1, p}_\lambda(\Omega)} 
\end{equation*} 
for $|\lambda| \ge \lambda_0$ and for $\lambda_0$ large enough. This completes the proof of  \eqref{thm-WP-2}. 

\medskip
\noindent {\underline{\it Step 2:}} Set 
\begin{multline*}
X = \Big \{(u_1,u_2) \in [L^p(\Omega)]^2  : \dive (A_1\nabla u_1),\dive (A_2\nabla u_2) \in L^p(\Omega), \\
u_1-u_2 \in W^{2,p}(\Omega), u_1-u_2 = 0 \mbox{ on } \Gamma, \mbox{ and }(A_1\nabla u_1 - A_2 \nabla u_2)\cdot \nu = 0 \mbox{ on } \Gamma \Big \}.
\end{multline*}
The space $X$ is a Banach space endowed with the norm
\begin{equation}
\|(u_1,u_2)\|_X : = \|(u_1,u_2)\|_{L^p(\Omega)} + \|\dive (A_1\nabla u_1),\dive (A_2\nabla u_2) \|_{L^p(\Omega)} + \|u_1-u_2\|_{W^{2,p}(\Omega)}.
\end{equation}

Define 
$$
B_\lambda : X \to [L^p(\Omega)]^2
$$  
by
\[
B_\lambda(u_1,u_2) = (\dive (A_1\nabla u_1)-\lambda \Sigma_1u_1, \dive (A_2\nabla u_2)-\lambda \Sigma_2u_2) .
\]
Clearly, $B_\lambda$ is bilinear and continuous on $X$. 

We claim that  
\begin{equation}\label{thm-WP-claim}
\mbox{$B_\lambda$ has a closed and dense range}. 
\end{equation}
Assuming this, we derive that 
\begin{equation}\label{thm-WP-claim-1}
B_\lambda(X) = [L^p(\Omega)]^2,  
\end{equation}
which yields the existence of the solutions. 

\medskip 
It remains to prove \eqref{thm-WP-claim}. 

\medskip 
We first prove that $B_\lambda$ has a closed range.  Let $((u_{1,n},u_{2,n}))_n \subset X$ be such that $(f_{1,n},f_{2,n}): = B_\lambda(u_{1,n},u_{2,n}) \to (f_1,f_2)$ in $[L^p(\Omega)]^2$. It follows from \eqref{thm-WP-1} by Step 1 that $((u_{1,n},u_{2,n}))_n$ is a Cauchy sequence in $X$. Let  $(u_1, u_2)$ denote  its limit. One can then show that $(f_{1,n},f_{2,n}) \to (f_1, f_2) : = B_\lambda(u_1,u_2) $  since $B_\lambda$ is continuous. Thus $B_\lambda$ has a closed range. 

\medskip 

We next establish that the range of $B_\lambda$ is dense.  To this end, it suffices to show that if   $(f_1,f_2) \in [L^q(\Omega)]^2$ with $\frac{1}{p}+\frac{1}{q} = 1$ is such that 
\begin{equation}
\label{thm-ex-111}
\int_{\Omega} \langle B_{\lambda}(u_1,u_2),(f_1,f_2)\rangle dx=0 \quad \quad \mbox{ for all }(u_1,u_2) \in X, 
\end{equation}
then $(f_1, f_2) = (0, 0)$.

Since \eqref{thm-ex-111} holds for all $(u_1, u_2) \in [C_c^\infty (\Omega)]^2 \subset X$, it follows that, for $\ell =1,2$,
\begin{equation}
\label{thm-ex-12}
\dive (A_\ell\nabla f_\ell) -\overline{\lambda}\Sigma_\ell f_\ell = 0 \mbox{ in }\Omega.
\end{equation}
Since $A_\ell \in C^1(\bar \Omega)$ and $f_\ell \in L^q(\Omega)$, using the standard regularity theory in $L^q$-scale, see also \cite[Lemma 17.1.5]{Hor3}, one has 
$$
f_\ell \in W^{2, q}_{\loc} (\Omega). 
$$

Set, in $\Omega$,  
\begin{equation}\label{thm-ex-def-g}
g_1 = f_1 \quad \mbox{ and } \quad g_2 = - f_2. 
\end{equation}
Then, by \eqref{thm-ex-12},  
\begin{equation}
\label{thm-ex-13}
\dive (A_\ell\nabla g_\ell) - \overline{\lambda}\Sigma_\ell g_\ell = 0 \mbox{ in } \Omega, 
\end{equation}
and, by \eqref{thm-ex-111},  for $(u_1,u_2) \in X$,
\begin{equation}\label{thm-ex-pp0}
\int_{\Omega} \dive(A_1 \nabla u_1) \bar g_1 - \lambda \Sigma_1 u_1 \bar g_1 -  \int_{\Omega} \dive(A_2 \nabla u_2) \bar g_2 - \lambda \Sigma_2 u_2 \bar g_2 = 0. 
\end{equation}

From \eqref{thm-ex-pp0}, we have, taking $(u_1,u_2) \in X \cap [W^{2,p}(\Omega)]^2$,
\begin{multline}\label{thm-ex-pp1}
\int_{\Omega} \dive\big(A_1 \nabla (u_1-u_2)  \big) \bar g_1 + \dive(A_1 \nabla u_2 ) (\bar g_1 - \bar g_2) + \dive  \big( (A_1- A_2) \nabla u_2 \big) \bar g_2
\\[6pt]- \lambda \Sigma_1 u_1 \bar g_1 + 
 \lambda \Sigma_2 u_2 \bar g_2  = 0.
\end{multline}
Using that $g_2 \in W^{2,q}_{loc}(\Omega)$ and $A_1 = A_2$ in $\Omega_\tau$, an integration by parts leads to
\begin{equation}\label{thm-ex-pp2}
\int_{\Omega} \dive  \big( (A_1- A_2) \nabla u_2 \big) \bar g_2 = \int_{\Omega}
\dive  \big( (A_1- A_2) \nabla \bar g_2 \big) u_2. 
\end{equation}
Since $u_1 - u_2 \in W^{2, p}(\Omega)$,  $u_1- u_2 = 0$ on $\Gamma$ and $A\nabla (u_1 - u_2) \cdot \nu = 0$ on $\Gamma$, there exists a sequence $(v_n)_n \subset C^2_c(\Omega)$ such that $v_n \to u_1 - u_2$ in $W^{2, p}(\Omega)$. An integration by parts yields 
\begin{multline}\label{thm-ex-pp3}
\int_{\Omega} \dive\big(A_1 \nabla (u_1-u_2)  \big) \bar g_1 
= \lim_{n \to + \infty} \int_{\Omega} \dive\big(A_1 \nabla v_n  \big) \bar g_1
\\[6pt]
= \lim_{n \to + \infty}  \int_{\Omega} \dive\big(A_1 \nabla \bar g_1  \big) v_n  =   \int_{\Omega} \dive\big(A_1 \nabla \bar g_1  \big) (u_1 - u_2). 
\end{multline}
Combining \eqref{thm-ex-pp1}, \eqref{thm-ex-pp2}, and \eqref{thm-ex-pp3} yields 
\begin{multline}\label{thm-ex-pp4}
\int_{\Omega}  \dive(A_1 \nabla u_2 ) (\bar g_1 - \bar g_2) 
= - \int_{\Omega}
\dive  \big( (A_1- A_2) \nabla \bar g_2 \big) u_2 - \int_{\Omega} \dive\big(A_1 \nabla \bar g_1  \big) (u_1 - u_2) \\[6pt]
+  \int_{\Omega} \lambda \Sigma_1 u_1 \bar g_1 -  \int_{\Omega}
 \lambda \Sigma_2 u_2 \bar g_2  \mathop{=}^{\eqref{thm-ex-13}} \int_{\Omega} \dive(A_1 \nabla (\bar g_1 - \bar g_2)) u_2. 
\end{multline}

Since $u_2$ can be chosen arbitrary\footnote{Taking then $u_1:= u_2$ so that $(u_1,u_2) \in X$.} in $W^{2, p}(\Omega)$,  and for every $\xi \in [C^1(\bar \Omega)]^d$ there exists $u_2 \in W^{2, p}(\Omega)$ with $u_2|_{\Gamma} =0$ such that $\dive(A_1 \nabla u_2 ) = \dive \xi $ with $\| u_2 \|_{W^{2, p}(\Omega)} \le C \|\xi \|_{L^p(\Omega)}$, 
it follows that, see,  e.g.,  \cite[Proposition 9.18]{Brezis-FA},  
$$
g_1 -  g_2  \in W^{1, q}_0 (\Omega). 
$$
This in turn implies, by \eqref{thm-ex-pp4},  that $g_1 -  g_2  \in W^{2, q}(\Omega)$ and $A \nabla ( g_1 - g_2) \cdot \nu = 0$ on $\Gamma$.  It follows that  $g_1 = g_2 = 0$ in $\Omega$ after applying  Step 1 to $(g_1, g_2)$ and  $\bar \lambda$ (instead of $\lambda$). Thus $f_1 = f_2 = 0$  by \eqref{thm-ex-def-g} and the proof of Step 2 is complete.  
\qed


\section{The Weyl law for the transmission eigenvalues} \label{sect-W}

Throughout this section, we assume  \eqref{condi3}. 
\subsection{The operator $T_\lambda$ and its adjoint $T_{\lambda}^*$.}

We first reformulate the Cauchy problem \eqref{ITE-C} in a form for which we can apply the theory of Hilbert-Schmidt operators. Given $(f, g) \in [L^p(\Omega)]^2$ ($1<p<+\infty$), assume that $(u_1,u_2) \in [L^p(\Omega)]^2$ with $u_1-u_2 \in W^{2,p}(\Omega)$ is a solution of \eqref{ITE-C}  with  $\lambda \in \C^*$ where, in $\Omega$, 
\[
f_1= \Sigma_1 f+\lambda^{-1}\Sigma_2 g\quad \quad \mbox{ and } \quad \quad f_2 = \lambda^{-1}\Sigma_2 g.
\]
Define, in $\Omega$, 
\[
u = u_1-u_2 \quad \quad \quad \mbox{ and }\quad \quad \quad v = \lambda u_2.
\]
Then the pair $(u,v) \in W^{2,p}(\Omega) \times L^p(\Omega)$ is a solution of 
\begin{equation}
\label{ITE-Cauchy-2}
\left \{
\begin{array}{rll}
\dive (A \nabla u) -\lambda \Sigma_1 u - (\Sigma_1-\Sigma_2)v &=   \Sigma_1 f  & \text{ in }\Omega ,\\[6pt]
\dive (A \nabla v) -\lambda \Sigma_2 v &=  \Sigma_2 g & \text{ in }\Omega ,\\[6pt]
u = 0, \: \: A\nabla u\cdot \nu &=0 & \text{ on }\Gamma.
\end{array}
\right .
\end{equation} 

As a direct consequence of \Cref{thm-WP} (see also \eqref{thm-WP-2}), we have

\begin{prop}
\label{prop-pre}
Assume \eqref{condi1}-\eqref{condi2}, and \eqref{condi3}-\eqref{condi5}. Let $c \in (0, 1)$ and $1<p<+\infty$. There exists $\lambda_0>0$ depending on $p$, $c$, $\Lambda$,  and $\Omega$ such that the following holds:  for $(f,g) \in [L^p(\Omega)]^2$
and for $\lambda \in \C$ with $|\Im (\lambda) | \geq c |\lambda| $ and $|\lambda|>\lambda_0$,  there exists a unique solution $(u,v) \in W^{2,p}(\Omega) \times L^p(\Omega)$ of the Cauchy problem \eqref{ITE-Cauchy-2};  moreover,  we have
\begin{equation}\label{prop-pre-p1}
 \|v\|_{L^p (\Omega) } +\|u\|_{W^{2,p}_{\lambda}(\Omega)} \leq  C|\lambda|^{-1/2} \left (    |\lambda|^{1/2}\|f\|_{L^p(\Omega)} + |\lambda|^{-1/2}\|g\|_{L^p(\Omega) } \right )
\end{equation}
and
\begin{equation}\label{prop-pre-p2}
\|v\|_{ W_{\lambda}^{1,p}(\Omega)} + \|u\|_{W^{3,p}_\lambda (\Omega)}\leq C  \left (  \|f\|_{W^{1,p}_\lambda(\Omega)}+|\lambda|^{-1/2}\| g \|_{ L^p(\Omega)} \right ), 
\end{equation}
for some positive constant $C$ independent of $\lambda$, $f$, and $g$. 
\end{prop}

As a consequence, we have 
\begin{cor}\label{cor-reg} Assume \eqref{condi1}-\eqref{condi2} and  \eqref{condi3}-\eqref{condi5}. Let  $c \in (0, 1)$, and $1 < p < + \infty$. There exists $\lambda_0>0$ depending on $p$, $c$, $\Lambda$,  and $\Omega$ such that the following holds:  for $(f, g)  \in W^{1, p}(\Omega) \times L^p(\Omega)$, and for $\lambda \in \C$ with  $|\Im (\lambda) | \geq c |\lambda| $ and $|\lambda|>\lambda_0$, there exists a unique solution $(u, v) \in W^{3,p}(\Omega) \times W^{1, p}(\Omega)$ of \eqref{ITE-Cauchy-2}; moreover, for 
\begin{enumerate}
\item either $1 < p < d$ and $p \le q \le \frac{dp}{d-p}$, 

\item either $d=p \le q < + \infty$, 

\item or $p> d$ and $q = + \infty$, 
\end{enumerate}
we have
\begin{equation}
\|v\|_{ L^q (\Omega)} +\|u\|_{W^{2,q}_\lambda(\Omega)}\leq  C |\lambda|^{\frac{d}{2} \left ( \frac{1}{p}-\frac{1}{q} \right ) -\frac{1}{2}} \left ( \|f\|_{W^{1,p}_\lambda (\Omega)}+|\lambda|^{-1/2} \|g\|_{ L^p(\Omega)} \right ), 
\end{equation}
for some positive constant $C$ independent of $\lambda$, $f$, and $g$. 
\end{cor}

\begin{rem}\label{rem-cont} \rm
In case (3) of \Cref{cor-reg}, we derive  that $(u,v) \in C^2(\overline{\Omega}) \times  C(\overline{\Omega})$.
\end{rem}

\begin{proof}
Choose $\lambda_0$ such that the conclusion of \Cref{prop-pre} holds. By Gagliardo-Nirenberg's interpolation inequalities (see \cite{Gagliardo59,Nirenberg59}), we have
\begin{equation*}
\| v\|_{L^q(\Omega)} \le C_{p, q, \Omega} \| v\|_{L^p(\Omega)}^{1 - a} \| v\|_{W^{1, p}(\Omega)}^{a} \leq C_{p, q, \Omega} \| v\|_{L^p(\Omega)}^{1 - a} \| v\|_{W^{1, p}_\lambda(\Omega)}^{a},  
\end{equation*}
where 
$$
a = d \left(\frac{1}{p} - \frac{1}{q} \right). 
$$
This implies 
\begin{equation*}
\| v\|_{L^q(\Omega)} \le C_{p, q, \Omega}   |\lambda|^{- \frac{1}{2} (1 - a)} \left (  \|f\|_{W^{1,p}_\lambda(\Omega)}+|\lambda|^{-1/2}\| g \|_{ L^p(\Omega)} \right ). 
\end{equation*}
The other assertions can be proved similarly.
\end{proof}

\begin{definition} Assume \eqref{condi1}-\eqref{condi2} and \eqref{condi3}-\eqref{condi5}.  Let $1 < p < + \infty$ and  $\lambda \in \mC$. System  \eqref{ITE-Cauchy-2} is said to be {\rm well-posed} in $L^p(\Omega) \times L^p(\Omega)$ if  the existence, the uniqueness, and \eqref{prop-pre-p1} and \eqref{prop-pre-p2} hold for $(f, g) \in L^p(\Omega) \times L^p(\Omega)$.  For $p=2$ and $\lambda \in \mC$ being such that   \eqref{ITE-Cauchy-2} is well-posed in $L^2(\Omega) \times L^2(\Omega)$, we define 
\begin{equation}
\label{pre-Tl}
\begin{array}{rccc}
T_\lambda : & L^2(\Omega) \times L^2(\Omega) &\to & L^2(\Omega) \times L^2(\Omega) \\[6pt]
& (f,g) & \mapsto & (u,v)
\end{array}
\end{equation}
where $(u,v)$ is the unique solution of \eqref{ITE-Cauchy-2}.
\end{definition}

\begin{rem}\rm Let $\lambda \in \mC$ satisfy the conclusion of \Cref{prop-pre} with $p=2$. Then  system \eqref{ITE-Cauchy-2} is well-posed in $L^2(\Omega) \times L^2(\Omega)$ and $T_\lambda$ is defined. 

\end{rem}

\begin{rem} \label{rem-comp}\rm Let $\lambda^* \in \C$ be such that $T_{\lambda}^*$ is defined.  If $\mu$ is a characteristic value of the operator $T_{\lambda^*}$ associated with an eigenfunction $(u,v)$ and if $\lambda^*+\mu \neq 0$ we have
\begin{equation}
\label{rem-comp-0}
\lambda^*+ \mu \mbox{ is a transmission eigenvalue of \eqref{ITE} } 
\end{equation} 
 with an eigenfunction pair $(u_1,u_2)$ given by 
\begin{equation}
\label{rem-comp-1}
u_1 = u + \frac{1}{\lambda^* + \mu} v \quad \quad \mbox{ and }\quad \quad u_2 = \frac{1}{\lambda^* + \mu} v.
\end{equation}
Moreover, the converse holds (see \Cref{rem-indep}). 
\end{rem}

\begin{rem} \rm Let $\lambda^* \in \C$ be such that $T_{\lambda}^*$ is defined. By \eqref{prop-pre-p1} and \eqref{prop-pre-p2} the range of $T_{\lambda^*}^2$ is a subset of $H^1(\Omega)\times H^1(\Omega)$.  It follows that the operator $T_{\lambda^*}^2$ is compact from $L^2(\Omega) \times L^2(\Omega)$ into itself. By the spectral theory of compact operators, see, e.g., \cite{Brezis-FA}, the spectrum of $T_{\lambda^*}^2$ consists in a discrete set of eigenvalues and the generalized eigenspace associated to each eigenvalue is of finite dimension.  As a consequence, the set of eigenvalues of $T_{\lambda^*}$ is discrete. This in turn implies that  the set of the transmission eigenvalues of \eqref{ITE} is discrete. This fact is previously established in \cite{MinhHung1} but the arguments presented here are different. 
\end{rem}

\begin{rem}\label{rem-indep}\rm 
Let $\lambda^* \in \C$ be such that $T_{\lambda}^*$ is defined.
If $\lambda_j$ is an eigenvalue of the transmission eigenvalue problem, then $\lambda_j \neq \lambda^*$ and $\lambda_j - \lambda^*$ is a characteristic value of $T_{\lambda^*}$.  One can show that the multiplicity of the characteristic values of $\lambda_j - \lambda^*$ and $\lambda_j - \hat{\lambda}$ associated with $T_{\lambda^*}$ and $T_{\hat{\lambda}}$ are the same. Hence the multiplicity of the eigenvalues associated with $T_{\lambda^*}$ is independent of $\lambda^*$. With this observation we define the multiplicity of $\lambda_j$ as the one of the characteristic value $\lambda_j - \lambda^*$ of $T_{\lambda^*}$. 
\end{rem}

The rest of this section is devoted to characterize the adjoint $T_\lambda^*$ of $T_\lambda$. This will be used in the proof of \Cref{prop-p}. To this end, for $(\widetilde{f},\widetilde{g}) \in [L^p(\Omega)]^2$ with $1 < p < + \infty$, we consider the system, for $(\widetilde{u}, \widetilde{v}) \in W^{1, p}(\Omega) \times L^p(\Omega)$   \footnote{We emphasize here that in the first equation of \eqref{ITE-Cauchy-3}, we have $\Sigma_2 \widetilde{u}$ not  $\Sigma_1 \widetilde{u}$, compare with \eqref{ITE-Cauchy-2}.}
\begin{equation}\label{ITE-Cauchy-3}
\left \{\begin{array}{cl}
\dive (A \nabla \widetilde{u}) -\lambda \Sigma_2 \widetilde{u} - (\Sigma_1-\Sigma_2)\widetilde{v} =   \Sigma_2 \widetilde{f}  & \text{ in }\Omega ,\\[6pt]
\dive (A \nabla \widetilde{v}) -\lambda \Sigma_1 \widetilde{v} =  \Sigma_1 \widetilde{g} & \text{ in }\Omega ,\\[6pt]
\widetilde{u} = 0, \quad  A\nabla \widetilde{u}\cdot \nu =0 & \text{ on }\Gamma.
\end{array}\right .
\end{equation}

Assume \eqref{condi1}-\eqref{condi2}, and \eqref{condi3}-\eqref{condi5}. Let $c \in (0, 1)$ and $1<p<+\infty$. By \Cref{prop-pre}, there exists $\lambda_0>0$ depending on $p$, $c$, $\Lambda$,  and $\Omega$ such that \eqref{ITE-Cauchy-3} is well-posed in $L^p(\Omega) \times L^p(\Omega)$ for $\lambda \in \C$ with $|\Im (\lambda) | \geq c |\lambda| $ and $|\lambda|>\lambda_0$, i.e., for $(f, g) \in L^p(\Omega) \times L^p(\Omega)$, there exists a unique solution $(\widetilde{u}, \widetilde{v}) \in W^{1,p}(\Omega) \times L^p(\Omega)$ of \eqref{ITE-Cauchy-3}; moreover,  
\begin{equation*}
\|\widetilde{u}\|_{W^{2,p}_{\lambda}(\Omega)}  + \|\widetilde{v}\|_{L^p (\Omega) }\leq  C|\lambda|^{-1/2} \left (    |\lambda|^{1/2}\|\widetilde{f}\|_{L^p(\Omega)} + |\lambda|^{-1/2}\|\widetilde{g}\|_{L^p(\Omega) } \right )
\end{equation*}
and
\begin{equation*}
\|\widetilde{u}\|_{W^{3,p}_\lambda (\Omega)} +  \|\widetilde{v}\|_{ W_{\lambda}^{1,p}(\Omega)} \leq C  \left (  \|\widetilde{f}\|_{W^{1,p}_\lambda(\Omega)}+|\lambda|^{-1/2}\| \widetilde{g} \|_{ L^p(\Omega)} \right). 
\end{equation*}

\begin{definition} Assume \eqref{condi1}-\eqref{condi2} and \eqref{condi3}-\eqref{condi5}.   For $p=2$ and for $\lambda \in \mC$ being such that   \eqref{ITE-Cauchy-3} is well-posed in $L^2(\Omega) \times L^2(\Omega)$, we define 
\begin{equation}
\begin{array}{rccc}
 \widetilde{T}_\lambda : & L^2(\Omega) \times L^2(\Omega) &\to & L^2(\Omega) \times L^2(\Omega) \\[6pt]
& (\widetilde{f},\widetilde{g})   & \mapsto & (\widetilde{u}, \widetilde{v})
\end{array}
\end{equation}
where $(\widetilde{u}, \widetilde{v})$ is the unique solution of \eqref{ITE-Cauchy-3}.
\end{definition}

\begin{lem}\label{lem-T*}
Assume \eqref{condi1}-\eqref{condi2} and \eqref{condi3}-\eqref{condi5}. Let $p=2$ and let $\lambda$ be such that $T_\lambda$ and $\widetilde{T}_{\bar \lambda}$ are defined. Set, for $x \in \Omega$,
\begin{equation}
\label{pre-7}
P(x) = 
\begin{pmatrix}
0 & \Sigma_1(x)  \\[6pt]
 \Sigma_2(x) & 0
\end{pmatrix}.
\end{equation} 
We have
\begin{equation}
T_\lambda^* = P \widetilde{T}_{\overline{\lambda}} P^{-1}.
\end{equation}
\end{lem}

\begin{proof} Fix $(f, g) \in [L^2(\Omega)]^2$ and $(f^*,g^*) \in [L^2(\Omega)]^2$.  Set $(u, v) = T_{\lambda}(f, g)$ and $(u^*,v^*) = \widetilde{T}_{\overline{\lambda}} P^{-1}(f^*,g^*)$.  Then 
\begin{equation}\label{lem-T*-p1}
\int_{\Omega} \langle (f,g), P \widetilde{T}_{\overline{\lambda}} P^{-1}(f^*,g^*) \rangle = \int_\Omega \Sigma_1 f \overline{v^*}+\Sigma_2 g \overline{u^*}. 
\end{equation}
Since $(u, v) = T_{\lambda}(f, g)$, we have 
\begin{equation}\label{lem-T*-p2}
\int_\Omega \Sigma_1 f \overline{v^*}+\Sigma_2 g \overline{u^*} 
= \int_{\Omega} (\dive (A\nabla u) -\lambda \Sigma_1 u - (\Sigma_1-\Sigma_2) v ) \overline{v^*} + \int_\Omega ( \dive (A\nabla v ) -\lambda \Sigma_2 v ) \overline{u^*}. 
\end{equation}
As in Step 2 of the proof of \Cref{thm-WP}, an integration by parts yields 
\begin{multline}\label{lem-T*-p3}
\int_{\Omega} (\dive (A\nabla u) -\lambda \Sigma_1 u - (\Sigma_1-\Sigma_2) v ) \overline{v^*} + \int_\Omega ( \dive (A\nabla v ) -\lambda \Sigma_2 v ) \overline{u^*} \\[6pt]
= \int_{\Omega} u (\overline{\dive (A\nabla v^*)-\overline{\lambda} \Sigma_1 v^*)} +   \int_{\Omega} v (\overline{\dive (A\nabla u^*)-\overline{\lambda} \Sigma_2 u^*-(\Sigma_1-\Sigma_2)v^*)}. 
\end{multline}
Since $(u^*,v^*) = \widetilde{T}_{\overline{\lambda}} P^{-1}(f^*,g^*)$, we have 
\begin{multline}\label{lem-T*-p4}
\int_{\Omega} u (\overline{\dive (A\nabla v^*)-\overline{\lambda} \Sigma_1 v^*)} +   \int_{\Omega} v (\overline{\dive (A\nabla u^*)-\overline{\lambda} \Sigma_2 u^*-(\Sigma_1-\Sigma_2)v^*)} \\[6pt]
=  \int_\Omega \langle T_\lambda (f,g), (f^*,g^*) \rangle .
\end{multline}
Combining \eqref{lem-T*-p1}-\eqref{lem-T*-p4} yields 
\begin{equation}
\int_{\Omega} \langle (f,g), P \widetilde{T}_{\overline{\lambda}} P^{-1}(f^*,g^*) \rangle =  \int_\Omega \langle T_\lambda (f,g), (f^*,g^*) \rangle, 
\end{equation}
and the conclusion follows. 
\end{proof}

\subsection{Hilbert-Schmidt operators} 

In this section, we recall the definition and several properties of Hilbert-Schmidt operators.  We begin with 
 
\begin{definition}
\label{def-HS}
Let $H$ be a separable Hilbert space and let $(\phi_j)_{j=1}^\infty$ be an orthonormal basis of $H$. 
\begin{enumerate}
\item Let $\cT$ be a linear and bounded operator on $H$. We say that $\cT$ is {\rm Hilbert-Schmidt} if its {\rm double norm} is finite, i.e.
\[
\vvvert \cT  \vvvert := \left ( \sum_{j=1}^\infty \|\cT \phi_j\|_H^2 \right )^{1/2} < + \infty.
\]
\item Let $\cT_1$ and $\cT_2$ be two Hilbert-Schmidt operators on $H$. The {\rm trace} of the composition $\cT_1 \cT_2$ is defined by
\[
{\rm trace}(\cT_1 \cT_2) := \sum_{j=1}^\infty ( \cT_1 \cT_2 \phi_j,\phi_j )_H.
\]
\end{enumerate}
\end{definition}

\begin{rem} \rm One can check that \Cref{def-HS} does not depend on the choice of the basis $(\phi_j)_{j=1}^\infty$ and the trace of $\cT_1 \cT_2$ is well defined as an absolutely convergent series (see \cite[Theorems 12.9 and  12.12]{Agmon}). 
\end{rem}

Let $m \in \mN$ and  $\bT : [L^2(\Omega)]^m \to [L^2(\Omega)]^m$ be a Hilbert-Schmidt operator. There exists a unique kernel  $\bK \in [L^2(\Omega \times \Omega)]^{m \times m}$, see e.g. \cite[Theorems 12.18 and  12.19]{Agmon},  such that 
\begin{align}\label{Z2}
	(	\mathbf{T}u)(x)=\int_\Omega \bK(x,y) u(y) dy ~~~\mbox{ for a.e. } \, x\in \Omega, \mbox{ for all } u \in [L^2(\Omega)]^m. 
	\end{align}
Moreover, 
\begin{equation}\label{norm-K-T}
\vvvert\mathbf{T} \vvvert^2=\mathop{\iint}_{\Omega\times\Omega}|\bK(x,y)|^2 \,  dx \,  dy. 
\end{equation}
Note that \cite[Theorems 2.18 and  12.19]{Agmon} state for  $m=1$,  nevertheless, the same arguments hold for $m \in \mN$ as noted in \cite{MinhHung2}. 

\medskip 
We have, see \cite{Agmon} (see also \cite[Lemma 4]{MinhHung2}):

\begin{lemma} \label{lem-T1T2} Let $m \in \mN$ and let 	$\mathbf{T}_1,\mathbf{T}_2$ be two Hilbert-Schmidt operators in $[L^2(\Omega)]^m$ with the corresponding kernels $\bK_1$ and $\bK_2$. Then $\bT:  = \mathbf{T}_1\mathbf{T}_2$ is a Hilbert-Schmidt operator with the kernel $\mathbf{K}$ given by 
\begin{equation}\label{lem-T1T2-1}
\bK(x, y) = \int_{\Omega} \bK_1 (x, z)  \bK_2 (z, y) \, dz. 
\end{equation}
Moreover,  
\begin{equation}\label{traceT^2}
	{\rm trace}  (\mathbf{T}_1\mathbf{T}_2)=\int_{\Omega} {\rm trace}(\mathbf{K}(x,x))dx.
	\end{equation}
\end{lemma}

We have, see, e.g., \cite[Lemma 3]{MinhHung2}. 
\begin{lemma}\label{lem-HS1} Let $d \ge 2$, $m \in \mN$,  and  $\bT : [L^2(\Omega)]^m \to [L^2(\Omega)]^m$ be  such that $\mathbf{T}(\phi)\in [C(\bar \Omega)]^m$ for $\varphi  \in [L^2(\Omega)]^m$,  and 
	\begin{equation}\label{lem-HS1-st0}
	\|\mathbf{T}(\phi) \|_{L^\infty (\Omega)}\leq M \|\phi \|_{L^2 (\Omega)}, 
\end{equation} 
for some $M \ge 0$. 
Then  $\mathbf{T}$ is a Hilbert-Schmidt operator,  
	\begin{align}\label{lem-HS1-st1}
\vvvert \mathbf{T} \vvvert \leq C_m |\Omega|^{1/2}M,
	\end{align}
	and  the kernel $\bK$ of $\bT$ satisfies
	\begin{equation}
	\label{lem-HS1-st1'}
	\sup_{x\in \Omega}\left(\int_{\Omega}|\bK(x,y)|^2  dy\right)^{1/2}\leq C_m |\Omega|^{1/2}M.
	\end{equation}
Assume in addition that 
		\begin{equation}\label{lem-HS1-st2}
	\| \mathbf{T}(\phi) \|_{L^\infty (\Omega)}\leq \widetilde{M}||\phi||_{L^1 (\Omega)} \mbox{ for } \phi\in [L^2(\Omega)]^m, 
	\end{equation}
	for some $\widetilde{M} \ge 0$, 
then the kernel $\bK$ of $\bT$ satisfies 
	\begin{align}\label{lem-HS1-st3}
	|\bK(x,y)|\leq  \widetilde{M} ~~ \mbox{ for a.e. } x, y \in \Omega.
	\end{align}
Here $C_m$ denotes a positive constant depending only on $m$. 
\end{lemma}

As a consequence of \Cref{lem-HS1}, we derive the following result. 

\begin{cor}
\label{cor-K-cont}
Let $\mathbf{T}_1$ and $\mathbf{T}_2$ two Hilbert-Schmidt operators on $[L^2(\Omega)]^m$ be such that the ranges of  $\mathbf{T}_1$ and $\mathbf{T}_2^*$ are in $[C(\bar{\Omega})]^m$ and \eqref{lem-HS1-st0} holds for $\mathbf{T}_1$ and $\mathbf{T}_2^*$. Assume that \eqref{lem-HS1-st2} holds for $\mathbf{T}= \mathbf{T}_1 \mathbf{T}_2$. Then the kernel $\mathbf{K}$ of $\mathbf{T}$ is continuous on $\bar \Omega \times \bar \Omega$ and \eqref{lem-HS1-st3} holds for every $(x,y) \in \Omega \times \Omega$.
\end{cor}

\begin{proof}
Let $\mathbf{K}_1$ (resp. $\mathbf{K}_2$) be  the kernel of $\mathbf{T}_1$ (resp. $\mathbf{T}_2$) and let $\mathbf{K}_2^*$ be the kernel of $\mathbf{T}_2^*$.  We claim that for  $\eps>0$,  there exists $\delta>0$ such that for every $(x, x')\in \Omega \times \Omega$ with $|x-x'|<\delta$ we have
\begin{equation}
\label{cor-K-cont-1}
 \left ( \int_{\Omega}|\mathbf{K}_1(x,z)-\mathbf{K}_1(x',z)|^2 dz \right )^{1/2} \le \eps \quad  \mbox{ and } \quad  \left ( \int_{\Omega}|\mathbf{K}_2^*(x,z)-\mathbf{K}_2^*(x',z)|^2 dz \right )^{1/2}  \le \eps. 
\end{equation}

 Admitting \eqref{cor-K-cont-1}, we continue the proof. 
We have, see, e.g., \cite[Theorem 12.20]{Agmon}, 
\begin{equation}
\label{cor-K-cont-11}
\mathbf{K}_2(z,y) = \overline{\mathbf{K}_2^*(y,z)}.
\end{equation}
Since
\[
\mathbf{K}(x,y) \mathop{=}^{\eqref{lem-T1T2-1},\eqref{cor-K-cont-11}}  \int_\Omega \mathbf{K}_1(x,z) \overline{\mathbf{K}_{2}^*(z,y)} dz, 
\]
it follows from \eqref{cor-K-cont-1} that ${\bf K}$ is continuous in $\bar \Omega \times \bar \Omega$. This in turn implies \eqref{lem-HS1-st3} by \Cref{lem-HS1} applied to $\mathbf{T}$.

It remains to prove \eqref{cor-K-cont-1}. We have
\begin{equation*}
\left ( \int_{\Omega}|\mathbf{K}_1(x,z)-\mathbf{K}_1(x',z)|^2 dz \right )^{1/2} \le C \sup_{\varphi \in [L^2(\Omega)]^m; \| \varphi\|_{L^2(\Omega)} \le 1} \left| \int_{\Omega}( \mathbf{K}_1(x,z)-\mathbf{K}_1(x',z) ) \varphi (z)  dz \right|. 
\end{equation*}
Given $\eps > 0$, let $\varphi_\eps \in [L^{2}(\Omega)]^m$ with $\|\varphi_\eps\|_{L^2(\Omega)}\leq 1$ be such that
\begin{equation*}
\left ( \int_{\Omega}|\mathbf{K}_1(x,z)-\mathbf{K}_1(x',z)|^2 dz \right )^{1/2}\leq \left | \int_\Omega (\mathbf{K}_1 (x, z)-\mathbf{K}_1(x', z))\varphi_\eps (z) dz \right | + \frac{\eps}{2}. 
\end{equation*}
This yields 
\begin{equation}
\label{cor-K-cont-3}
\left ( \int_{\Omega}|\mathbf{K}_1(x,z)-\mathbf{K}_1(x',z)|^2 dz \right )^{1/2}\leq  |(\mathbf{T}_1 \varphi_\eps)(x)-(\mathbf{T}_1 \varphi_\eps)(x')|  + \frac{\eps}{2}. 
\end{equation}
The first inequality of \eqref{cor-K-cont-1} now follows from \eqref{cor-K-cont-3} and the fact that ${\bf T}_1 \varphi_\eps \in [C(\bar \Omega)]^m$. 

\medskip 
Similarly, we obtain the second inequality of \eqref{cor-K-cont-1}.  
\end{proof}

\subsection{The operators  $\bT_{\theta, t}$ and their properties}
Denote 
\begin{equation}
\label{def-k}
k= \left [\frac{d}{2} \right ] + 1, 
\end{equation}
the smallest integer greater than $d/2$. Fix
\begin{equation}
2=p_1 < p_2 < \cdots < p_{k}<+\infty
\end{equation}
such that 
\begin{equation}
p_{j-1}<p_j < \frac{dp_{j-1}}{d-p_{j-1}} \quad  \mbox{ and } \quad p_{k}>d.
\end{equation}

Denote 
\begin{equation}
\label{def-lambda^*}
\lambda^* = t^*e^{i\frac{\pi}{2}}, 
\end{equation}
for some large $t^*>0$ such that,  for $t \ge t^*$,  \eqref{ITE-Cauchy-2} with $\lambda =  t e^{i\frac{\pi}{2}}$ is well-posed in $L^p(\Omega) \times L^p(\Omega)$ and  \eqref{ITE-Cauchy-3} with $\lambda = t e^{- i\frac{\pi}{2}}$ is well-posed in $L^p(\Omega) \times L^p(\Omega)$
 with $p =p_1, \, \cdots, \, p_{k}$.

Let 
\begin{equation}\label{omega-j}
\mbox{$\omega_j \in \C$ with $1 \le j \le k+1$ be the (distinct) $(k+1)$-th roots of $1$ (thus $\omega_j^{k+1} = 1$)}
\end{equation} 
and let
\begin{equation}
\label{def-gam}
\Theta = \R \setminus \left \{ \frac{\pi}{k+1}\mathbb{Z} \right \}.
\end{equation}

\begin{definition}\label{def-t-theta}
For $\theta \in \Theta$, $1 \le j \le k+1$,  and $t>0$,  we define
\begin{equation}
\label{p4}
\lambda_{j, \theta , t} = \lambda^* +  \omega_j  t e^{i\theta},
\end{equation}
and \[t_\theta>t^*\] such that the following properties hold, for $t \ge t_\theta$,  
\begin{multline}
\label{unif-t0-1}
\mbox{ \eqref{ITE-Cauchy-2} with $\lambda =  \lambda_{j, \theta , t} $ is well-posed in $L^p(\Omega) \times L^p(\Omega)$} \\
\mbox{ and  \eqref{ITE-Cauchy-3} with $\lambda = \bar \lambda_{j, \theta , t} $ is well-posed in $L^p(\Omega) \times L^p(\Omega)$ with $p =p_1, \, \cdots, \, p_{k}$}, 
\end{multline}
and
\begin{equation}
\label{unif-t0-2}
\frac{t}{2} \leq |\lambda_{j, \theta, t}| < 2t. 
\end{equation}
\end{definition}

Such a $t_\theta>t^*$ exists by \Cref{prop-pre} after noting that, for  $\theta \in \Theta$, 
\[
\Im \left (\omega_j e^{i\theta} \right ) \neq 0,
\]
and, for  $1 \le j \le k+1$,  
\begin{equation*}
\lim_{t \to + \infty} \frac{\left |\Im \left (\lambda^*+ t\omega_j e^{i\theta} \right )\right |}{\left |\lambda^*+ t \omega_j e^{i\theta} \right |} = \left |\Im \left (\omega_j e^{i\theta} \right ) \right |>0.
\end{equation*}


Viewing \eqref{prop-pre-p1} and \eqref{prop-pre-p2}, it is convenient to modify $T_{\lambda_{j, \theta, t}}$ to capture the scaling with respect to $t \sim \lambda_{j,\theta,t}$ there, as in \cite{Robbiano16}.  Denote 
\begin{equation}\label{def-Mt}
M_t = 
\begin{pmatrix}
t^{1/2} & 0 \\
0 & t^{-1/2}
\end{pmatrix}.
\end{equation}
Let $\theta\in \Theta$   and $t \ge t_\theta$. Define,  for $1 \le j \le k+1$, 
\begin{equation}
\label{p11} 
T_{j, \theta, t} = M_t T_{\lambda_{j, \theta, t}} M_{t}^{-1} \quad \quad  \mbox{ and } \quad \quad \bT_{\theta, t} = T_{k+1, \theta, t}\circ T_{k, \theta, t} \circ \cdots \circ T_{1, \theta, t}.
\end{equation}

Here is the main result of this section.

\begin{prop} 
\label{prop-p} Let $\theta\in \Theta$ and let $t_\theta$ be given in \Cref{def-t-theta}. Then, for $t \ge t_\theta$,  
\begin{equation}
\label{prop-p-0}
\|\bT_{\theta,t}\|_{L^2(\Omega)\to L^2(\Omega)} \leq Ct^{-k-1}, 
\end{equation}
the range of $\bT_{\theta, t}$ is in $[C(\bar{\Omega})]^2$,  
\begin{equation}
\label{prop-p-4}
\|\bT_{\theta, t}\|_{L^2 (\Omega) \to L^\infty (\Omega)} \leq C t^{-k-1 + \frac{d}{4}}, 
\end{equation}
and
\begin{equation}
\label{prop-p-4*}
\|\bT_{\theta, t}\|_{L^1 (\Omega) \to L^2 (\Omega)} \leq C t^{-k-1+\frac{d}{4}}
\end{equation}
for some positive constant $C$ independent of $t$. Similar facts  hold for $\bT_{\theta, t}^*$. 
\end{prop}

As a direct consequence of \Cref{lem-HS1} and \Cref{prop-p},  we obtain 
\begin{cor}\label{cor-p} Let $\theta\in \Theta$   and 
 let $t_\theta$ be given in \Cref{def-t-theta}. Then, for $t \ge t_\theta$, the operator $\bT_{\theta, t}$ is Hilbert-Schmidt,   and 
\begin{equation}
\label{prop-p-5}
\vvvert  \bT_{\theta, t} \vvvert \leq Ct ^{-k-1+ \frac{d}{4}}, 
\end{equation} 
for some positive constant $C$ independent of $t$. 
\end{cor}

We now give 

\begin{proof}[Proof of \Cref{prop-p}] We first deal with 
\eqref{prop-p-4}. 
By using \eqref{unif-t0-2}, we derive that 
\[
\|T_{j,\theta,t}\|_{L^2(\Omega) \to L^2(\Omega)}\mathop{\leq}^{{\rm \Cref{prop-pre}}} Ct^{-1}
\]
and hence
\[
\|\bT_{\theta,t}\|_{L^2(\Omega) \to L^2(\Omega)} \leq \prod_{j=1}^{k+1}  \|T_{j,\theta,t}\|_{L^2(\Omega) \to L^2(\Omega)} \leq Ct^{-k-1}.
\]
This establishes \eqref{prop-p-0}.

Next we deal with \eqref{prop-p-4}. For $j=1,\cdots,k+1$ and $(f,g) \in [L^2(\Omega)]^2$, we write
\begin{equation}
\label{p12}
(u^{(j)},v^{(j)}) = T_{j, \theta, t}\circ T_{j-1, \theta,t} \circ \cdots \circ T_{1, \theta, t} (f,g).
\end{equation}
By \eqref{unif-t0-2}, we have
\begin{equation}
\label{prop-p-1}
t^{-1/2}\|u^{(1)}\|_{W^{1,2}_t(\Omega)} + \| v^{(1)} \|_{L^2(\Omega)} \mathop{\leq}^{{\rm \Cref{prop-pre}}} C t^{-1}  \|(f, g)\|_{L^2(\Omega) \times L^2(\Omega)}, 
\end{equation}
and,  for $2 \le j \le k$, 
\begin{multline}
\label{prop-p-2}
t^{-1/2}\|u^{(j)}\|_{W^{1, p_j}_t(\Omega)}+\|v^{(j)}\|_{ L^{p_j}(\Omega)} \\[6pt]
 \mathop{\leq}^{{\rm\Cref{cor-reg}}} C t^{-1+\frac{d}{2} \left ( \frac{1}{p_{j-1}} -\frac{1}{p_j} \right )}
\Big(  t^{-1/2} \|u^{(j-1)}\|_{W^{1,p_{j-1}}_t(\Omega)} + \|v^{(j-1)}\|_{L^{p_{j-1}}(\Omega)} \Big), 
\end{multline}
and 
\begin{multline}
\label{prop-p-3}
t^{-1/2}\|u^{(k+1)}\|_{W^{1,\infty}_t(\Omega)} +\|v^{(k+1)} \|_{L^\infty(\Omega)} \\[6pt] 
 \mathop{\leq}^{{\rm\Cref{cor-reg}}}  C t^{-1 + \frac{d}{2p_k}}
\Big(  t^{-1/2} \| u^{(k)}\|_{W^{1,p_{k}}_t(\Omega)} + \| v^{(k)}\|_{L^{p_{k}}(\Omega)} \Big).
\end{multline}
We derive from \eqref{prop-p-1}, \eqref{prop-p-2} and \eqref{prop-p-3},
\begin{multline}
\label{prop-p-6}
\| u^{(k+1)}\|_{L^\infty(\Omega)}+\|v^{(k+1)}\|_{L^\infty(\Omega) }  \\[6pt]
 \le C  t^{-1}  t^{-1 + \frac{d}{2p_k}} \prod_{j=2}^k t^{-1+\frac{d}{2} \left ( \frac{1}{p_{j-1}}-\frac{1}{p_j} \right )} \|(f,g)\|_{L^2(\Omega) \times L^2(\Omega)} = C t^{-k-1 + \frac{d}{4}} \|(f,g)\|_{L^2(\Omega) \times L^2(\Omega)}.  
\end{multline}
Thus \eqref{prop-p-4} is proved. 

We next establish \eqref{prop-p-4*}. We have, by \Cref{lem-T*}, 
\[
T_{j, \theta, t}^* =  M_t^{-1} P \widetilde{T}_{\overline{\lambda}_{j, \theta,t}} P^{-1} M_t, 
\]
where $P$ is given by \eqref{pre-7}. This implies 
\[
\bT_{\theta, t}^* = M_{t}^{-1} P \widetilde{T}_{\overline{\lambda}_{1, \theta,t}} \circ  \cdots  \circ \widetilde{T}_{\overline{\lambda}_{k+1, \theta, t}} P^{-1}M_t. 
\]

Similarly to \eqref{prop-p-6}, we have
\begin{equation}
\label{prop-p-7}
\|\bT_{\theta, t}^*\|_{L^2 (\Omega) \to L^\infty(\Omega)} \leq C t^{-k-1 + \frac{d}{4}}.
\end{equation}
By a standard dual argument,  we  derive from \eqref{prop-p-7} that 
\[
\|\bT_{\theta, t}\|_{L^1(\Omega) \to L^2 (\Omega) } \leq C t^{-k-1 + \frac{d}{4}}.
\]

The properties for ${\bf T}_{\theta, t}$ are established. 

\medskip 
The properties for ${\bf T}_{\theta, t}^*$ can be derived similarly. 
\end{proof}

\subsection{The approximation of the trace of a kernel}  

Denote 
\begin{equation}\label{p1}
\alpha = \frac{\pi}{4(k+1)} \quad \quad \mbox{ and } \quad \quad \beta = \frac{5\pi}{4(k+1)}.
\end{equation}
Then 
\begin{equation}\label{p1*}
\alpha,\beta \in \Theta \quad \mbox{ and }\quad  e^{i\alpha(k+1)}+e^{i\beta(k+1)} = 0.
\end{equation}
Recall that $\Theta$ is defined in \eqref{def-gam}. 



\begin{lem}  
\label{lem-Tt}
For $t \ge \max\{t_\alpha, t_\beta \}$,  where $t_\alpha$ and $t_\beta$ are given in \Cref{def-t-theta}, we have
\begin{enumerate}

\item the operator $\bT_{\alpha,t} \bT_{\beta,t}$ is Hilbert-Schmidt, and
\begin{equation}
\label{lem-Tt-1}
\vvvert  \bT_{\alpha,t} \bT_{\beta,t} \vvvert \leq C t^{-2k-2+\frac{d}{2}}; 
\end{equation}

\item the range of $\bT_{\alpha,t} \bT_{\beta,t}$ is in $[C(\bar{\Omega})]^2$, and 
\begin{equation}
\label{lem-Tt-2}
\|\bT_{\alpha,t} \bT_{\beta,t}\|_{L^1(\Omega) \to L^\infty(\Omega)} \leq Ct^{-2k-2+\frac{d}{2}}; 
\end{equation}

\item the kernel $\bK_t$ of $\bT_{\alpha,t} \bT_{\beta,t}$ is continuous in $\Omega \times \Omega$,  and 
\begin{equation}
\label{lem-Tt-3}
|\bK_t(x,y)| \leq C t^{-2k-2+\frac{d}{2}} \quad \mbox{ for all }(x,y) \in \Omega \times \Omega; 
\end{equation}
\end{enumerate}
for some positive constant $C$ independent of $t$. 
\end{lem}

\begin{proof}
Assertion \eqref{lem-Tt-1} follows from  \Cref{cor-p} and 
\[
\vvvert  \bT_{\alpha,t} \bT_{\beta,t} \vvvert \leq \vvvert  \bT_{\alpha , t} \vvvert \vvvert  \bT_{\beta, t} \vvvert.
\]
Applying \Cref{prop-p} and using the fact
\[
\|\bT_{\alpha,t} \bT_{\beta,t}\|_{L^1(\Omega) \to L^\infty(\Omega)} \leq \|\T_{\alpha,t} \|_{L^2(\Omega) \to L^\infty(\Omega)}\|\T_{\beta,t} \|_{L^1(\Omega) \to L^2(\Omega)}, 
\]
we obtain \eqref{lem-Tt-2}. 

Since both the range of $\bT_{\alpha,t}$ and $\bT_{\beta,t}^*$ are contained in $[C(\bar{\Omega})]^2$, the continuity of $\bK_t$ and \eqref{lem-Tt-3} follow from \Cref{cor-K-cont} and \eqref{lem-Tt-2}.
\end{proof}
 
For $\ell =1, \, 2$, $\theta \in \Theta$,  and $t>1$, consider, with $\lambda = t e^{i\theta}$, 
\begin{equation}
\label{def-S}
\begin{array}{rccc}
S_{\ell,\lambda,x_0} :&L^2(\R^d)& \to &L^2(\R^d)  \\[6 pt]
& f_\ell & \mapsto & v_\ell
\end{array}
\end{equation}
where $v_\ell \in H^1(\R^d)$ is  the unique solution of
\begin{equation}\label{vl-Rd}
\dive (A(x_0)\nabla v_\ell) -\lambda \Sigma_\ell (x_0) v_\ell = \Sigma_\ell (x_0) f_\ell \mbox{ in } \mR^d. 
\end{equation}
One then has
\begin{equation}
S_{\ell,\lambda,x_0} f(x) = \int_{\R^d}F_{\ell,\lambda}(x_0,x-y)f(y) dy, 
\end{equation}
where
\begin{equation}
\label{def-F}
F_{\ell,\lambda}(x_0,z)= -\frac{1}{(2\pi)^d}\int_{\R^d} \frac{e^{i z \cdot \xi}}{\Sigma_{\ell}(x_0)^{-1}A(x_0) \xi \cdot \xi+\lambda}d\xi.
\end{equation}

Set,  for $\ell =1,2$, 
\[
\mathcal{S}_{\ell,t,x_0} = S_{\ell,\lambda_{k+1,\alpha,t},x_0}\circ \cdots \circ S_{\ell,\lambda_{1,\alpha,t},x_0} \circ S_{\ell,\lambda_{k+1,\beta,t},x_0}\circ \cdots \circ S_{\ell,\lambda_{1,\beta,t},x_0}.
\] 
Define,  for $\ell =1,2$, 
\begin{equation}
\label{ker-33}
\mathcal{F}_{\ell, t}(x_0,z) = \frac{1}{(2\pi)^d} \int_{\R^d} \frac{e^{i z \xi} \, d \xi}{\prod_{j=1}^{k+1} \left ( \Sigma_\ell (x_0)^{-1} A(x_0) \xi \cdot \xi + \lambda_{j,\alpha, t} \right ) \left ( \Sigma_\ell (x_0)^{-1} A(x_0) \xi \cdot \xi + \lambda_{j,\beta, t} \right )}. 
\end{equation}
Then 
\begin{equation}
\label{def-mS}
\mathcal{S}_{\ell,t,x_0} f_\ell (x) = \int_{\R^d}\mathcal{F}_{\ell,t}(x_0,x-y)f_\ell (y) dy.
\end{equation}
Since $2k+2>d$, the integrand appearing in \eqref{ker-33} belongs to $L^1(\R^d) \cap L^2(\R^d)$, and thus
\begin{equation}
z \mapsto \mathcal{F}_{\ell, t}(x_0,z) \mbox{ is continuous and belongs to }L^2(\R^d).
\end{equation}

To introduce the freezing coefficient version of \eqref{ITE-Cauchy-2} in the whole space, we use the following result in which \eqref{lem-rd-S} is the system of $(u, v) : = (v_1 - v_2, \lambda v_1)$, where $v_\ell$ ($\ell = 1, 2$) is defined by \eqref{vl-Rd}. 

\begin{lem}
 \label{lem-rd}
 Let $x_0 \in \Omega$, $c\in (0,1)$, $\lambda \in \C$ with $|\lambda|\ge 1$ and $|\Im (\lambda) |\geq c|\lambda|$. Let $p>1$ and let $(f,g) \in [L^p(\R^d)]^2$. Then there exists a unique solution  $(u,v) \in [W^{1,p}(\R^d)]^2$  of
 \begin{equation}
 \label{lem-rd-S}
 \left \{
 \begin{array}{rll}
   \dive (A(x_0)\nabla u) -\lambda \Sigma_1(x_0)u-(\Sigma_1(x_0)-\Sigma_2(x_0))v&=  \Sigma_1 (x_0) f & \mbox{ in } \R^d, \\ [6 pt]
      \dive ( A(x_0)\nabla v) -\lambda \Sigma_{2}(x_0) v &= \Sigma_2 (x_0) g & \mbox{ in }\R^d.
 \end{array}
 \right .
 \end{equation}
Moreover, 
\begin{equation}
 \label{lem-rd-E}
 \|u\|_{W^{2,p}_\lambda(\R^d)}+|\lambda|^{-1}\|v\|_{W^{2,p}_\lambda (\R^d)} 
 \leq C \left ( \|f\|_{L^p(\R^d)} +|\lambda|^{-1}\|g\|_{L^p(\R^d)} \right ),
 \end{equation}
 for some $C>0$ depending only on $\Lambda$,$c$ and $p$. As a consequence, 
\begin{enumerate}
\item either $1 < p < d$ and $p \le q \le \frac{dp}{d-p}$, 

\item either $d=p \le q < + \infty$, 

\item or $p> d$ and $q = + \infty$, 
\end{enumerate}
we have 
\[
 \|u\|_{W^{1,q}_\lambda (\R^d)} + |\lambda|^{-1}\|v\|_{W^{1,q}_\lambda (\R^d)} \leq  C |\lambda|^{\frac{d}{2} \left ( \frac{1}{p}-\frac{1}{q} \right ) -\frac{1}{2}}  \left ( \|f\|_{L^p(\R^d)} + |\lambda|^{-1}\|g\|_{L^p(\R^d)} \right ).  
\]
\end{lem}

\begin{proof} We emphasize here that \eqref{lem-rd-S} is a system with constant coefficients imposed in $\mR^d$.  The proof is quite standard. The idea is first to obtain the existence, uniqueness, and the estimate for $v$ using the second equation of \eqref{lem-rd-S},  and then using these to derive the ones for $u$ using the first equation of \eqref{lem-rd-S}. The details are omitted. 
\end{proof}

For $x_0 \in \Omega$, $j=1,\cdots,k+1$,  $\theta \in \Theta$, and $t>1$,  define, for $1 < p < + \infty$, 
\[
\begin{array}{cccc}
R_{\lambda_{j, \theta, t},x_0} :& [L^p(\R^d)]^2 &  \to &  [L^p(\R^d)]^2 \\[6 pt]
& (f, g) & \mapsto &   (u, v)
\end{array}
\]
where $(u, v) \in [W^{1,p}(\R^d)]^2$ is the unique solution \eqref{lem-rd-S} with $\lambda = \lambda_{j, \theta,t}$. Recall that $\lambda_{j, \theta,t}$ is defined in \eqref{p4}. We also introduce
\begin{equation}
\label{ker-8}
R_{j, \theta , t,x_0} = M_t R_{\lambda_{j, \theta,t}, x_0}M_t^{-1} \quad \quad \mbox{ and }\quad \quad \bR_{\theta, t, x_0} = R_{k+1, \theta, t, x_0}\circ \cdots \circ R_{1, \theta, t, x_0}. 
\end{equation}

As in the proof of \Cref{prop-p},  however, using \Cref{lem-rd}  instead of
\Cref{prop-pre} and \Cref{cor-reg}, we obtain 

\begin{lem}\label{lem-proS} Let $\theta \in \Theta$ and $t>1$.  Then,  the range of $\bR_{\theta,t,x_0}$ and $\bR_{\theta,t,x_0}^*$ are in $[C(\R^d)]^2$ for all $t>1$. Moreover, 
\begin{equation}
\label{ker-37}
\| \bR_{\theta,t,x_0}\|_{L^2 (\mR^d) \to L^\infty (\mR^d) }+ \|\bR_{\theta,t,x_0}\|_{L^1 (\mR^d) \to L^2 (\mR^d)} \leq Ct^{-k-1 + \frac{d}{4}}
\end{equation}
and
\begin{equation}
\label{ker-371}
\| \bR_{\theta,t,x_0}^*\|_{L^2 (\mR^d) \to L^\infty (\mR^d) }+ \|\bR_{\theta,t,x_0}^*\|_{L^1 (\mR^d) \to L^2 (\mR^d)} \leq Ct^{-k-1 + \frac{d}{4}}, 
\end{equation}
for some positive constant $C$ independent of $t$. 
\end{lem}

Define
\begin{equation}
\bR_{t,x_0} = \bR_{\alpha,t,x_0}\bR_{\beta,t,x_0}. 
\end{equation}
One can then write $\bR_{t,x_0} $ under the form 
\[
\bR_{t,x_0} 
=
\begin{pmatrix}
(\bR_{t,x_0})_{11} & (\bR_{t,x_0})_{12} \\
(\bR_{t,x_0})_{21} & (\bR_{t,x_0})_{22} \\
\end{pmatrix}.
\]
Note that, by the definition of $S_{\ell, \lambda, x_0}$, 
\[
R_{\lambda_{j,\theta,t},x_0}
= 
\begin{pmatrix}
S_{1,\lambda_{j,\theta,t},x_0} & \: \Sigma_1(x_0)^{-1}(\Sigma_1(x_0)-\Sigma_2(x_0))S_{1,\lambda_{j,\theta,t},x_0} S_{2,\lambda_{j,\theta,t},x_0} \\
0&\: S_{2,\lambda_{j,\theta,t},x_0}
\end{pmatrix}.
\]
It follows that  $R_{\lambda_{j,\theta,t},x_0}$ is an upper triangular matrix operator,  and  so is $\bR_{t,x_0}$. We deduce that
\[
(\bR_{t,x_0})_{21} = 0
\]
and,  for $\ell=1,2$, 
\[
(\bR_{t,x_0})_{\ell \ell} = \mathcal{S}_{\ell,t,x_0}.  
\]
These simple observations are useful in computing the approximation of the trace of the kernel of $\bR_{t, x_0}$.

As an immediate consequence of \eqref{def-mS}, $\bR_{t,x_0}$ is an integral operator whose  kernel verifies, for $\ell = 1, 2$,  
\begin{equation}\label{kernel-Kll}
(\bK_{t,x_0})_{\ell \ell}(x,y)=\mathcal{F}_{\ell}(x_0,x-y) \mbox{ for } x, y \in \mR^d. 
\end{equation}

Further properties of $\bK_{t,x_0}$ are given in  the following lemma.
\begin{lem} 
\label{lem-ker0} Let $t \ge 1$ and $x_0 \in \Omega$. Then $\bK_{t,x_0}$ is continuous on $\R^d \times \R^d$, and, for  $(x,y) \in \R^d\times \R^d$, it holds, for $\ell =1, 2$,  
\begin{equation}\label{lem-ker0-st1}
|(\bK_{t,x_0})_{\ell \ell}(x,y)| \leq Ct^{-2k-2+\frac{d}{2}}.
\end{equation}
Moreover,
\begin{multline}\label{lem-ker0-st2}
\mbox{trace}(\bK_{t,x_0} (x_0,x_0)) \\[6pt]
= \frac{t^{-2k-2+\frac{d}{2}} }{(2\pi)^d} \sum_{\ell = 1}^2 \int_{\R^d} \frac{d\xi}{(\Sigma_\ell (x_0)^{-1}A(x_0)\xi \cdot \xi)^{2k+2}-i}  + o(t^{-2k-2+\frac{d}{2}}) \mbox{ as } t \to + \infty.
\end{multline}
\end{lem}

\begin{remark} \rm Assertion \eqref{lem-ker0-st2} holds uniformly with respect to $x_0 \in \Omega$. 
\end{remark}

\begin{proof}  From \eqref{kernel-Kll}, it follows that  $(\bK_{t,x_0})_{\ell \ell}(x,y)$ is continuous on $\R^d\times \R^d$. By the choice of $\alpha$, $\beta$, and $\omega_j$ in \eqref{p1}, \eqref{p1*},  and \eqref{omega-j},  one has 
\begin{multline*}
\prod_{j=1}^{k+1} \left ( \Sigma_\ell (x_0)^{-1} A(x_0) \xi \cdot \xi + \lambda^* + \omega_j t e^{i\alpha}  \right ) \left ( \Sigma_\ell (x_0)^{-1} A(x_0) \xi \cdot \xi +  \lambda^* + \omega_j t e^{i\beta}    \right ) \\[6pt]
= (\Sigma_\ell(x_0)^{-1} A(x_0) \xi \cdot \xi + \lambda^*)^{2(k+1)} - i t^{2(k+1)}. 
\end{multline*}
It follows from \eqref{ker-33} that,    for every $x_0 \in \Omega$ and every $z \in \R^d$, 
\[
\mathcal{F}_{\ell,t}(x_0,z) = \frac{1}{(2\pi)^d} \int_{\R^d} \frac{e^{i z\cdot \xi} \, d \xi}{ (\Sigma_\ell(x_0)^{-1} A(x_0) \xi \cdot \xi + \lambda^*)^{2(k+1)} - i t^{2(k+1)}}. 
\]
A change of variables yields 
\begin{equation}\label{ker-41}
\mathcal{F}_{\ell,t}(x_0,z) = \frac{t^{-2k-2+\frac{d}{2} }}{(2\pi)^d} \int_{\R^d} \frac{e^{i t^{1/2}z\cdot  \xi} \, d \xi}{ (\Sigma_\ell(x_0)^{-1} A(x_0) \xi \cdot \xi + t^{-1}\lambda^*)^{2(k+1)} - i}. 
\end{equation}
Assertion \eqref{lem-ker0-st1} follows from \eqref{ker-41} since $|e^{i t^{1/2}z\cdot  \xi}|=1$ and $\lambda^* t^{-1}$ is uniformly bounded with respect to $t\geq 1$. 

By taking $z = 0$ in \eqref{ker-41}, we obtain \eqref{lem-ker0-st2} after using the dominated convergence theorem. 

\medskip 
The proof is complete. 
\end{proof}

We now prove the main result of this section concerning the trace of ${\bf T}_{\alpha, t} {\bf T}_{\beta, t}$ where 
$\alpha, \beta$ are given in \eqref{p1} and  ${\bf T}_{\theta, t}$ is defined in \eqref{p11}.

\begin{prop}
\label{prop-ker}
We have
\[
{\rm trace} (\bT_{\alpha, t} \bT_{\beta, t} ) = {\bf c} t^{-2k-2+ \frac{d}{2}} + o (t^{-2k-2+ \frac{d}{2}}) \quad \mbox{ as } \quad t \to + \infty,
\]
where
\begin{equation}
\label{prop-ker-2}
{\bf c} = \frac{1}{(2\pi)^d} \sum_{\ell=1}^2 \int_\Omega \int_{\R^d} \frac{d\xi \, dx}{\left ( \Sigma_\ell^{-1}(x) A(x) \xi \cdot \xi \right )^{2k+2}-i} . 
\end{equation}
\end{prop}

The proof of \Cref{prop-ker} uses the following result. 

\begin{lem}
\label{lem-approx} Let $\delta_0 \in (0,1)$ and $\theta \in \Theta$.  For  every $\eps > 0$, there exists $\delta_\eps \in (0,\delta_0/2)$ depending on $\eps$ such that the following holds: There exists $t_{\eps}>0$ depending on $\eps$ and $\delta_\eps$  such that for every $t > t_{\eps}$ and every $x_0 \in \Omega \setminus \overline{\Omega_{\delta_0}}$,  we have
\begin{equation}
\label{lem-approx-p1}
\|\bT_{\theta, t} -\bR_{\theta,t, x_0} \mathds{1}_\Omega \|_{L^2(\Omega) \to L^\infty(B(x_0,\delta_\eps))} \leq  \epsilon  t^{-k-1 + \frac{d}{4}}
\end{equation}
and
\begin{equation}\label{lem-approx-p2}
\|\bT_{\theta, t} \mathds{1}_\Omega  -\bR_{\theta, t, x_0}  \|_{L^2(\R^d) \to L^\infty(B(x_0,\delta_\eps))} \leq \eps  t^{-k-1 + \frac{d}{4}}, 
\end{equation}
and similar facts for $\bT_{\theta, t}^*$ and $\bR_{\theta, t, x_0}^*$.
\end{lem}

Recall that ${\bf R}_{\theta, t, x_0}$ is defined in \eqref{ker-8}.  We admit \Cref{lem-approx} and give the proof of 
\Cref{prop-ker}. The proof of \Cref{lem-approx} is presented right after the one of \Cref{prop-ker}.

\begin{proof}[Proof of \Cref{prop-ker}]
For $\epsilon>0$, let $\delta_0>0$ be such that
\begin{equation}
\label{ker-1}
|\Omega_{2\delta_0}|<\epsilon, 
\end{equation}
where $\Omega_\tau$ is given in \eqref{def-Omegatau}. 

We claim  that  there exists $\tau_*>0$,  depending  on $\Omega$,  and $\eps$ but independent of $x_0$, and a positive constant $C$, independent of $\eps$ and $x_0$, such that,  for $t> \tau_*$, 
\begin{equation}
\label{ker-42}
|\mbox{trace} (\bK_t (x_0, x_0))- \mbox{trace}({\bf K}_{t,x_0}(x_0,x_0))| \leq C \epsilon t^{-2k-2+\frac{d}{2}} \quad  \mbox{ for } x_0 \in \Omega \setminus \Omega_{\delta_0}. 
\end{equation}

Indeed, let $\chi \in C_c^\infty(\mR^d)$ be such that $\chi = 1$ in $B_1$ and $\supp \chi \subset B_2$. Denote, for $\delta \in (0,\delta_0/10)$,  
$$
\chi_{\delta, x_0} = \chi \big( ( \cdot - x_0) / \delta \big),
$$   
and define 
\begin{equation}
\label{ker-441}
\left\{\begin{array}{c}
{\rm \bold{P}}_{1,t,\delta} = \chi_{\delta, x_0}   (\bR_{\alpha,t,x_0} \mathds{1}_\Omega -\bT_{\alpha,t}) \bT_{\beta,t} \chi_{\delta, x_0}, \\[6pt]
{\rm \bold{P}}_{2,t,\delta} =  \chi_{\delta, x_0} \bR_{\alpha,t,x_0}  (\bR_{\beta,t,x_0} - \mathds{1}_\Omega \bT_{\beta,t} )\chi_{\delta, x_0}. 
\end{array} \right. 
\end{equation}

Then 
\begin{equation}
\label{ker-43}
\chi_{\delta, x_0} (\bR_{\alpha,t,x_0} \bR_{\beta,t,x_0} -\bT_{\alpha,t} \bT_{\beta,t}  ) \chi_{\delta, x_0} = {\rm \bold{P}}_{1,t,\delta} + {\rm \bold{P}}_{2,t,\delta}.  
\end{equation}

 By applying \Cref{lem-approx} below  with $\theta \in \{\alpha,\beta\}$, there exist $\delta_\eps>0$ and $t_{\eps}>0$ depending on $\eps$  such that for every $t>t_{\eps}$, 
 \begin{equation}
\label{ker-421}
\|  \chi_{\delta_\eps, x_0}   (\bT_{\alpha, t}-\bR_{\alpha, t, x_0} \mathds{1}_\Omega )  \|_{L^2(\Omega) \to L^\infty(\Omega)} \leq \epsilon  t^{-k-1 + \frac{d}{4}}
\end{equation}
and 
\begin{equation}
\label{ker-422}
\|  \chi_{\delta_\eps, x_0}  (   \bT_{\beta, t}^* \mathds{1}_\Omega-\bR_{\beta, t, x_0}^* )   \|_{L^2(\R^d) \to L^\infty(\Omega)} \leq \epsilon  t^{-k-1 + \frac{d}{4}}. 
\end{equation}
Since 
\[
\Big( ( \mathds{1}_\Omega \bT_{\beta,t} -\bR_{\beta, t, x_0} ) \chi_{\delta_\eps, x_0}  \Big)^* =  \chi_{\delta_\eps, x_0}  (   \bT_{\beta, t}^* \mathds{1}_\Omega-\bR_{\beta, t, x_0}^*),   
\]
  we derive from \eqref{ker-422}, using a dual argument, that 
\begin{equation}
\label{ker-48}
\| (\mathds{1}_\Omega \bT_{\beta,t}-\bR_{\beta,t,x_0}  ) \chi_{\delta_\eps, x_0}  \|_{L^1(\Omega) \to L^2(\mR^d)} \le \epsilon t^{-k-1+\frac{d}{4}}.
\end{equation}

By \Cref{prop-p} and \Cref{lem-proS}, we have 
\begin{equation}
\label{ker-481}
\| \bT_{\beta,t} \chi_{\delta_\eps, x_0} \|_{L^1(\Omega) \to L^2 (\Omega)} +  \| \chi_{\delta_\eps, x_0} \bR_{\alpha,t,x_0} \|_{L^2 (\mR^d) \to L^\infty (\Omega)}   \le C t^{-k-1 + \frac{d}{4}}
\end{equation}
for some constant $C>0$ independent of $\epsilon$ and $t$.

Using the fact, for appropriate linear operators $L_1$ and $L_2$,  
$$
\| L_1 L_2 \|_{L^1(\Omega) \to L^\infty(\Omega)} \le \| L_1 \|_{L^2 (\Omega) \to L^\infty (\Omega)} \| L_2 \|_{L^1(\Omega) \to L^2 (\Omega)},  
$$
and
$$
\| L_1 L_2 \|_{L^1(\Omega) \to L^\infty(\Omega)} \le \| L_1 \|_{L^2 (\R^d) \to L^\infty (\Omega)} \| L_2 \|_{L^1(\Omega) \to L^2 (\R^d)},  
$$
we derive from \eqref{ker-421}, \eqref{ker-48}, and \eqref{ker-481} that 
\begin{equation*}
\| {\rm \bold{P}}_{1,t,\delta_\eps}  + {\rm \bold{P}}_{2,t,\delta_\eps} \|_{L^1(\Omega)\to L^\infty(\Omega)}  \le 
\| {\rm \bold{P}}_{1,t,\delta_\eps}\|_{L^1(\Omega)\to L^\infty(\Omega)}  + \|{\rm \bold{P}}_{2,t,\delta_\eps} \|_{L^1(\Omega)\to L^\infty(\Omega)} \leq  C \epsilon t^{-2k-2+\frac{d}{2}}. 
\end{equation*}
This yields, by  \eqref{ker-43}, 
\begin{equation}
\label{ker-49}
\| \chi_{\delta_\eps, x_0}   (\bR_{\alpha,t,x_0} \bR_{\beta,t,x_0} -\bT_{\alpha,t} \bT_{\beta,t}  ) \chi_{\delta_\eps, x_0}  \|_{L^1(\Omega)\to L^\infty(\Omega)} \leq  C \epsilon t^{-2k-2+\frac{d}{2}}. 
\end{equation}

Since, for $x \in \Omega$, $\ell =1,2$ and $f\in L^2(\Omega)$, 
\begin{multline*}
\chi_{\delta_\eps, x_0}  \Big((\bR_{\alpha,t,x_0} \bR_{\beta,t,x_0})_{\ell \ell} -(\bT_{\alpha,t} \bT_{\beta,t})_{\ell \ell}  \Big) \chi_{\delta_\eps, x_0}  f (x) \\
= \chi_{\delta_\eps, x_0}  (x) \int_\Omega \chi_{\delta_\eps, x_0}  (y) \Big( (\bK_{t,x_0})_{\ell \ell}(x,y)- (\bK_t)_{\ell \ell} (x,y) \Big)f(y) dy,   
\end{multline*}
it follows that $\chi_{\delta_\eps, x_0}  (x) \chi_{\delta_\eps, x_0}  (y) ( (\bK_{t,x_0})_{\ell \ell}(x,y)- (\bK_t(x,y))_{\ell \ell})$ is the kernel of the operator 
$$
\chi_{\delta_\eps, x_0}   \Big((\bR_{\alpha,t,x_0} \bR_{\beta,t,x_0})_{\ell \ell} -(\bT_{\alpha,t} \bT_{\beta,t})_{\ell \ell}  \Big) \chi_{\delta_\eps, x_0}.
$$ 
By \Cref{lem-Tt} and \Cref{lem-ker0}, this kernel is continuous  on $\bar \Omega \times \bar \Omega$. 
Using \eqref{ker-49} and applying \Cref{lem-HS1}, we derive that, since   $\chi_{\delta_\eps, x_0}  (x_0) = 1$, 
\begin{equation}
\label{ker-50}
|\mbox{trace}(\bK_{t}(x_0,x_0))- \mbox{trace}(\bK_{t, x_0}(x_0,x_0))| \leq C \epsilon t^{-2k-2+ \frac{d}{2}} \quad  \mbox{ for all }t>t_{\eps}.
\end{equation}
Since the LHS of \eqref{ker-50} does not depend on $\eps>0$, the claim \eqref{ker-42} is proved. 

By \Cref{lem-Tt} we have, for $t>0$ large enough,
\begin{equation}
\label{ker-501}
\int_{\Omega_{2\delta_0}} |\mbox{trace}(\bK_t(x,x))| dx \leq C |\Omega_{2\delta_0} |t^{-2k-2+ \frac{d}{2}} \mathop{\leq}^{\eqref{ker-1}} C\epsilon t^{-2k-2+ \frac{d}{2}}
\end{equation}
and, similarly by \Cref{lem-ker0},
\begin{equation}\label{ker-502}
\int_{\Omega_{2\delta_0}} |\mbox{trace}(\bK_{t, x}(x, x))| dx \leq  C\epsilon t^{-2k-2+ \frac{d}{2}}.
\end{equation}

Combining  \eqref{ker-42},\eqref{ker-501},  and \eqref{ker-502} yields 
\begin{equation} \label{ker-521}
 \int_\Omega |\mbox{trace}(\bK_t(x,x))- \mbox{trace}(\bK_{t, x}(x, x))| dx  \leq C\eps t^{-2k-2+\frac{d}{2}} \quad \mbox{ for all }t>t_\eps.
\end{equation}
The conclusion follows from \Cref{lem-ker0} and \eqref{ker-521}.
\end{proof}

We now give 

\begin{proof}[Proof of \Cref{lem-approx}]  Let $\eps >0 $ and $\theta \in \Theta$. 
First, we prove \eqref{lem-approx-p1}. Fix $\chi \in C_c^\infty (\R^d)$ such that $\supp \chi \subset B_2$ and  $\chi = 1$  in $B_1$. Set, for $0< \delta < \delta_0/ 100$, 
\[
\chi_\delta = \chi \big( (\cdot -x_0)/\delta \big) \mbox{ in } \mR^d. 
\]

Define, for  $(f,g) \in [L^2(\Omega)]^2$, and for $j=1,\cdots, k+1 = [d/2] + 2$, 
\[
(u^j,v^j) = T_{\lambda_{j, \theta,t}}\circ \cdots \circ T_{\lambda_{1, \theta, t}}  (u^0, v^0)  \quad  \mbox{ and } \quad (u^j_0,v^j_0) = S_{\lambda_{j, \theta, t},x_0}\circ \cdots \circ S_{\lambda_{1, \theta, t},x_0}  (u^0_0,v^0_0), 
\]
where 
\begin{equation}
\label{case1}
(u^0,v^0)=(f, g) \quad \quad \mbox{ and } \quad \quad (u^0_0,v^0_0)=  (\mathds{1}_\Omega f, \mathds{1}_\Omega g) \quad \mbox{ in } \Omega. 
\end{equation}
Set, for $0 \le j \le k+1$,  
$$
(u^{j, \delta}, v^{j, \delta}) = (\chi_\delta u^j, \chi_\delta v^j)
\quad \quad \mbox{ and } \quad \quad 
(u^{j, \delta}_0, v^{j, \delta}_0) = (\chi_\delta u^j_0, \chi_\delta v^j_0). 
$$


We have
\begin{equation}
\label{iter}
(u^{j, \delta}, v^{j, \delta}) = S_{\lambda_{j, \theta, t}, x_0} (u^{j-1, \delta}, v^{j-1, \delta}) + S_{\lambda_{j, \theta, t}, x_0} (f^{j, \delta}, g^{j, \delta}), 
\end{equation}
where
\begin{multline}
\label{cc3}
\Sigma_1 (x_0) f^{j, \delta} =  (\Sigma_1 (x) -\Sigma_1(x_0)) u^{j-1, \delta} -\lambda_{j, \theta,t} (\Sigma_1(x_0)-\Sigma_1(x)) u^{j, \delta}+ 
A(x)\nabla \chi_\delta \cdot \nabla u^j \\[6pt]  - (\Sigma_1(x_0)-\Sigma_1(x)- \Sigma_2(x_0)+\Sigma_2(x))  v^{j, \delta} + \dive \Big( (A(x_0)-A(x)) \nabla u^{j, \delta} + u^j A(x)\nabla \chi_\delta\Big)
\end{multline}
and 
\begin{multline}
\label{cc4}
\Sigma_2(x_0) g^{j, \delta}  = (\Sigma_2 (x) -\Sigma_2(x_0)) v^{j-1, \delta} - \lambda_{j, \theta,t} (\Sigma_2(x_0)-\Sigma_2(x)) v^{j, \delta} + A(x) \nabla \chi_\delta \cdot \nabla v^j \\[6pt]
+ \dive \Big( (A(x_0)-A(x)) \nabla v^{j, \delta} + v^j A(x)\nabla \chi_\delta  \Big). 
\end{multline}

Similarly, we have
\[
(u^{j, \delta}_0, v^{j, \delta}_0) = S_{\lambda_{j, \theta, t}, x_0} (u^{j-1, \delta}_0, v^{j-1, \delta}_0) + S_{\lambda_{j, \theta, t},x_0} (f^{j, \delta}_0, g^{j, \delta}_0), 
\]
where 
\[
\Sigma_1(x_0) f^{j, \delta}_0 =
A(x_0) \nabla \chi_\delta \cdot \nabla u^j_0 + \dive \Big(u^j_0 A (x_0)\nabla \chi_\delta\Big)
\]
and 
\[
\Sigma_2(x_0) g^{j, \delta}_0  =  A(x_0) \nabla \chi_\delta \cdot \nabla v^j_0 \\[6pt]
+ \dive \Big(  v^j_0 A(x_0) \nabla \chi_\delta  \Big). 
\]

For $r>0$, define
\[
\Phi (r) = \min \left \{1, \sup_{|x-y|<r} \left( |A(x)-A(y)|+ \sum_{\ell=1}^2  |\Sigma_\ell(x)-\Sigma_\ell(y)| \right) \right \}.
\] 
We claim that 
\begin{multline}
\label{cc1}
\|f^{j, \delta}\|_{L^{p_{j-1}}(\Omega \setminus \Omega_{\delta_0/2})} + t^{-1}\|g^{j, \delta}\|_{L^{p_{j-1}}(\Omega \setminus \Omega_{\delta_0/2})} \\[6pt] \le  C_{\delta_0}  \left(\Phi(\delta)  +  \frac{1}{\delta t^{1/2}} + \frac{1}{\delta^2 t} \right) \left ( \|u^{j-1}\|_{L^{p_{j-1}}(\Omega)} + t^{-1}\|v^{j-1}\|_{L^{p_{j-1}}(\Omega)} \right ), 
\end{multline}
and
\begin{multline}
\label{cc2}
\|f^{j, \delta}_0\|_{L^{p_{j-1}}(\Omega \setminus \Omega_{\delta_0/2})} + t^{-1}\|g^{j, \delta}_0\|_{L^{p_{j-1}}(\Omega \setminus \Omega_{\delta_0/2})} \\[6pt] \le  C_{\delta_0}  \left( \frac{1}{\delta t^{1/2}} + \frac{1}{\delta^2 t} \right) \left ( \|u^{j-1}_0\|_{L^{p_{j-1}}(\Omega)} + t^{-1}\|v^{j-1}_0\|_{L^{p_{j-1}}(\Omega)} \right ). 
\end{multline}

We first admit  \eqref{cc1} and \eqref{cc2} and continue the proof.
Since, in $\Omega$,
$$
(u^{0, \delta}, v^{0, \delta}) = (u^{0, \delta}_0, v^{0, \delta}_0),  
$$
using \eqref{iter}, \eqref{cc1} and \eqref{cc2} and \Cref{lem-rd},  for $j=1$ and then for $j=2, \dots, k+1$,  we have 
\begin{multline}\label{lem-approx-coucou}
\|u^{j, \delta} - u^{j, \delta}_0\|_{L^{p_{j}}(\Omega \setminus \Omega_{\delta_0/2})} + t^{-1}\|v^{j, \delta} - v^{j, \delta}_0\|_{L^{p_{j}}(\Omega \setminus \Omega_{\delta_0/2})} \\[6pt]
\le C_{\delta_0} \left(\Phi (\delta)  +  \frac{1}{\delta t^{1/2}} + \frac{1}{\delta^2 t} \right) t^{- \frac{d}{2 p_j}- j + \frac{d}{4}  } \left(  \| f\|_{L^2(\Omega)} + t^{-1 } \| g\|_{L^2(\Omega)} \right).
\end{multline}
Fix  $\delta = \delta_\eps>0$ such that $C_{\delta_0} \Phi(\delta_\eps)<\eps/2$. Take  $t_\eps>0$ sufficiently large  such that $C_{\delta_0}(\delta_\eps^{-1}t^{-1/2} + \delta^{-2}_\eps t^{-1})<\eps/2$ for every $t>t_\eps$. Taking $j = k+1$ in \eqref{lem-approx-coucou} gives \eqref{lem-approx-p1}. 

The proof of  \eqref{lem-approx-p2} is similar to the one of  \eqref{lem-approx-p1} by considering 
$(u^0, v^0)$ and $(u^0_0, v^0_0)$ defined as follows  
\[
(u^0,v^0)= \mathds{1}_\Omega  (f, g) \quad \quad \mbox{ and } \quad \quad (u^0_0,v^0_0)= (f, g) \quad \mbox{ in } \mR^d, 
\]
instead of \eqref{case1}. 

\medskip 
Similar facts for $\bT_{\theta, t}^*$ and $\bR_{\theta, t, x_0}^*$ by analogous arguments. 

\medskip 

It remains to establish \eqref{cc1} and \eqref{cc2}.   
From the definition of $(u^j, v^j)$ and the theory of elliptic equations (see e.g. \cite[Theorem 9.11]{GilbargTrudinger}),  we have, for $\Omega_1 \Subset \Omega_2 \subset \Omega$,   
\begin{equation}
\label{fact1}
\|u^{j}\|_{W^{2,p}_t(\Omega_1)} +t^{-1}\|v^{j}\|_{W^{2,p}_t (\Omega_1)}   \leq C \left ( \|u^{j-1}\|_{L^p(\Omega_2)} + t^{-1}\|v^{j-1}\|_{L^p(\Omega_2)} \right )
\end{equation}
and, similarly,
\begin{equation}
\label{fact2}
\|u^{j}_0\|_{W^{2,p}_t(\Omega_1)} + t^{-1}\|v^{j}_0\|_{W^{2,p}_t (\Omega_1)}   \leq C \left ( \|u^{j-1}_0\|_{L^p(\Omega_2)} + t^{-1}\|v^{j-1}_0\|_{L^p(\Omega_2)} \right ), 
\end{equation}
for some positive constant $C$ independent of $f$, $g$, and $t$. It follows that 
\begin{equation}
\label{fact3}
 \|\nabla u^{j, \delta}\|_{L^p(\Omega_1)} +t^{-1}\|\nabla v^{j, \delta}\|_{L^p (\Omega_1)}  \leq   C \left( \frac{1}{\delta t} +  \frac{1}{t^{1/2}} \right)  \left ( \|u^{j-1}\|_{L^p(\Omega_2)} + t^{-1}\|v^{j-1}\|_{L^p(\Omega_2)} \right ),
\end{equation}
\begin{multline}
\label{fact4}
 \|\nabla^2 u^{j, \delta}\|_{L^p(\Omega_1)} +t^{-1}\|\nabla^2 v^{j, \delta}\|_{L^p (\Omega_1)}  \\[6pt]
  \leq   C \left( 1 + \frac{1}{\delta t^{1/2}} +  \frac{1}{\delta^2 t} \right)  \left ( \|u^{j-1}\|_{L^p(\Omega_2)} + t^{-1}\|v^{j-1}\|_{L^p(\Omega_2)} \right ),
\end{multline}
and
\begin{equation}
\label{fact5}
 \|u^{j, \delta}_0\|_{W^{1, p}_t(\Omega_1)} +t^{-1}\|v^{j, \delta}_0\|_{W^{1,p}_t (\Omega_1)}  \\[6pt] 
 \leq    C \left( \frac{1}{\delta t} +  \frac{1}{t^{1/2}} \right)   \left ( \|u^{j-1}_0\|_{L^p(\Omega_2)} + t^{-1}\|v^{j-1}_0\|_{L^p(\Omega_2)} \right ).
\end{equation}

By \eqref{cc3} and \eqref{cc4}, we have
\begin{align*}
C \left (  \right. & \left.  \|f^{j,\delta}\|_{L^{p_{j-1}}(\Omega_1)}+ t^{-1}\|g^{j,\delta}\|_{L^{p_{j-1}}(\Omega_1)} \right ) \le  \Phi (\delta) \Big( \|u^{j-1}\|_{L^{p_{j-1}}}+t^{-1}\|v^{j-1}\|_{L^{p_{j-1}}(\Omega_1)} \Big )  \\
& + \Big (  \Phi(\delta) + \frac{1}{\delta^2 t}+\frac{1}{\delta t^{1/2}}\Big ) \Big ( \|u^j\|_{W^{2,p_{j-1}}_t(\Omega_1)}+t^{-1} \|v^j\|_{W^{2,p_{j-1}}_t(\Omega_1)} \Big ) \\
& +\|\nabla u^{j,\delta}\|_{L^{p_{j-1}}(\Omega_1)}+ t^{-1}\|\nabla v^{j,\delta}\|_{L^{p_{j-1}}(\Omega_1)}  +\Phi(\delta) \Big ( \| \nabla^2 u^{j,\delta}\|_{L^{p_{j-1}}(\Omega_1)} +  t^{-1}\|\nabla^2 v^{j,\delta}\|_{L^{p_{j-1}}(\Omega_1)} \Big ). 
\end{align*}
Combining \eqref{fact1}-\eqref{fact5} yields \eqref{cc1}. Estimate \eqref{cc2} follows  similarly. 

\medskip 
The proof is complete. 
\end{proof}

\subsection{A connection of the counting function and the trace of $\bT_{\alpha, t} \bT_{\beta, t}$ for large $t$}

We start this section by recalling the definition of the modified resolvent of an operator (see, e.g., \cite[Definition 12.3]{Agmon}).
\begin{definition} \label{def-modRes}
Let $H$ be a Hilbert space and $\mathcal{T} : H \to H$ be a linear and bounded operator. The modified resolvent set $\rho_m(\mathcal{T})$ of $\mathcal{T}$ is the set of all non-zero complex numbers $s$ such that $I-s  \mathcal{T}$ is bijective and $(I-s \mathcal{T})^{-1}$ is bounded on $H$. For $s \in \rho_m(\mathcal{T})$ the transformation $(\mathcal{T})_s = \mathcal{T}(I-s\mathcal{T})^{-1}$ is the modified resolvent of $\mathcal{T}$.
\end{definition}

Recall that, for $s \in \rho_m(\mathcal{T})$, we have 
\begin{equation}\label{pro-mod}
(\mathcal{T})_s = \mathcal{T}(I-s\mathcal{T})^{-1} = (I-s\mathcal{T})^{-1}\mathcal{T}.
\end{equation}

Let $\mathcal{T}: L^2(\Omega) \times L^2(\Omega) \to L^2(\Omega) \times L^2(\Omega)$ be a linear and bounded operator. We have for $z \in \C$ (see e.g. \cite{Robbiano13})
\begin{equation}
\label{HS-3}
I-z^{k+1}\mathcal{T}^{k+1} = \prod_{j=1}^{k+1} (I-\omega_jz\mathcal{T})
\end{equation}
and 
\begin{equation}
\label{HS-4}
I-z^{k+1}\mathcal{T}^{k+1} \mbox{ is invertible } \Longleftrightarrow   I-\omega_jz\mathcal{T} \mbox{ is invertible for every }j.
\end{equation}
Recall that $\omega_j^{k+1}=1$. Using the decomposition \eqref{HS-3},  and the equivalence in \eqref{HS-4}, one can prove the following lemma.
\begin{lem}
\label{lem-mod-res} 
Let $\widetilde{\theta} \in \R\setminus \{\pi \mathbb{Z}\}$. Set $\theta := \frac{\tilde{\theta}}{k+1} \in \Theta$. 
There exists $t_\theta > 1$ such that,  for every  $t>t_\theta$,  it holds 
\begin{equation}
\label{lem-mod-res-2}
\gamma:=t^{k+1} e^{i\widetilde{\theta}} \in \rho_m(T_{\lambda^*}^{k+1})
\end{equation}
and 
\begin{equation}
\label{lem-mod-res-3}
(T_{\lambda^*}^{k+1})_\gamma = M_t^{-1} \mathbf{T}_{\theta,t} M_t. 
\end{equation}
\end{lem}

\begin{proof} We have, by the definition of $\gamma$, 
\begin{equation}
\label{lem-mod-res-21}
I- \gamma T_{\lambda^*}^{k+1} \mathop{=}^{\eqref{HS-3}} \prod_{j=1}^{k+1} (I-\omega_j t e^{i\theta} T_{\lambda^*}).
\end{equation}
By \Cref{prop-pre}, there exists $t_\theta > 0$ such that $T_{\lambda_* + \omega t e^{i \theta}}$ is defined for $t \ge t_\theta$.  Hence 
\begin{equation}
\label{lem-mod-res-211}
\omega_j t e^{i\theta} \in \rho_m(T_{\lambda^*}) \quad \mbox{ and } \quad  (T_{\lambda^*})_{\omega_j t e^{i\theta}} = T_{\lambda^* + \omega_j t e^{i\theta}} = T_{\lambda_{j,\theta,t}} \mbox{ for }t\geq t_\theta.
\end{equation}
(see, e.g. \cite[Lemma 3.1]{Fornerod-Ng1} for the arguments in a similar setting).
Combining \eqref{lem-mod-res-21} and \eqref{lem-mod-res-211} leads  $\gamma \in \rho_m (T_{\lambda^*}^{k+1})$ for $t \ge t_\theta$.

It follows from \eqref{pro-mod} that, for $t \ge t_\theta$, 
\begin{equation}
\label{HS-6}
T_{\lambda_{j,\theta,t}} = T_{\lambda^*} (I-\omega_j t e^{i\theta} T_{\lambda^*})^{-1}=(I-\omega_j t e^{i\theta} T_{\lambda^*})^{-1}T_{\lambda^*}
\end{equation}
and thus,
\begin{multline}
M_t^{-1} \mathbf{T}_{\theta,t}M_t \mathop{=}^{\eqref{p11}} \prod_{j=1}^{k+1} T_{\lambda_{j,\theta,t}} \mathop{=}^{\eqref{HS-6}} T_{\lambda^*}^{k+1}\prod_{j=1}^{k+1} (I-\omega_j e^{i\theta}t T_{\lambda^*})^{-1}  \\
=  T_{\lambda^*}^{k+1} (I-\gamma T^{k+1}_{\lambda^*})^{-1} \mathop{=}^{\rm \mbox{def.}} (T_{\lambda^*}^{k+1})_{\gamma}.
\end{multline}
The proof is complete.
\end{proof}

The following proposition establishes a connection between the trace of the operator $\bT_{\alpha,t}\bT_{\beta,t}$ and the counting function for large  $t$.
The arguments of the proof are in the spirit of \cite{MinhHung2} (see also \cite{Robbiano13}).

\begin{prop} 
\label{prop-spec} We have
\[
\mathcal{N}(t) = \frac{\Im ({\bf c})}{\frac{d}{8(k+1)} \int_0^\infty s^{\frac{d}{8(k+1)}-1}(1+s)^{-1}ds} t^{\frac{d}{2}} + o(t^{\frac{d}{2}})\quad \quad  \mbox{ as }t \to + \infty, 
\]
where $\mathbf{c}$ is given by \eqref{prop-ker-2}. 
\end{prop}

\begin{proof} 
%
%

For $t$ sufficiently large, by \Cref{lem-mod-res}, we have 
$$
(T_{\lambda^*}^{k+1})_{t^{k+1} e^{ i (k+1) \alpha}} = M_t^{-1} \mathbf{T}_{\alpha, t} M_t. 
$$
Note that  
\begin{equation}\label{HS-Rm}
\Big((T_{\lambda^*}^{k+1})_{\gamma_1} \Big)_{\gamma_2} = (T_{\lambda^*}^{k+1})_{\gamma_1 + \gamma_2} 
\end{equation}
provided that $\gamma_1$,$\gamma_1 + \gamma_2 \in \rho_m(T_{\lambda^*}^{k+1})$. 
It follows from \Cref{lem-mod-res} that, for large $t$ and for $s \ge 0$,  
$$
-2 (t+s)^{k+1} e^{i (k+1) \alpha} \in \rho_m \Big( M_t^{-1} \mathbf{T}_{\alpha, t} M_t \Big). 
$$

Let $s_{1}, s_{2}, \dots$ be the characteristic values of $M_t^{-1} \mathbf{T}_{\theta,t} M_t$  repeated a number of times equal to their multiplicities. Applying \cite[Theorem 12.17]{Agmon}, we have 
\begin{multline}
\label{HS-1500}
\mbox{trace}\Big(M_t^{-1}\bT_{\alpha,t} M_t (M_t^{-1}\bT_{\alpha,t}M_t )_{-2 (t+s)^{k+1} e^{i (k+1) \alpha}} \Big) = \sum_{j} \frac{1}{s_{j}(s_{j}+2e^{i\alpha (k+1)} (t+s)^{k+1})}+c_t. 
\end{multline}

We claim that 
\begin{equation}\label{prop-spec-ct}
c_t = 0. 
\end{equation}

Assume this, we continue the proof.  As a consequence of \eqref{HS-1500} with $s=0$, we have 
\begin{equation}
\label{HS-1500*}
\mbox{trace}\Big(M_t^{-1}\bT_{\alpha,t} M_t (M_t^{-1}\bT_{\alpha,t}M_t )_{-2t^{k+1} e^{i (k+1) \alpha}} \Big) = \sum_{j} \frac{1}{s_{j}(s_{j}+2e^{i\alpha (k+1)} t^{k+1})}. 
\end{equation}

Let $(\mu_j)_j$ be the set of characteristic values of $T_{\lambda^*}$ repeated according to their multiplicity. It is well-known that $\mu_j^{k+1}$ are the characteristic values of $T_{\lambda^*}^{k+1}$ and the multiplicity of $\mu_j^{k+1}$ is equal to the sum of the one of the characteristic values $\mu$ of  $T_{\lambda^*}$  such that $\mu^{k+1} = \mu_j^{k+1}$. By \Cref{lem-mod-res}, for large $t$, $e^{i \alpha (k+1)}t^{k+1}$ is not a characteristic value of $T_{\lambda^*}^{k+1}$. We obtain, by \cite[Theorem 12.4]{Agmon}, that the set of the characteristic values of  $(T_{\lambda^*}^{k+1})_{t^{k+1} e^{ i (k+1) \alpha}}$ is given by 
$$
\Big\{ \mu_j^{k+1} - t^{k+1} e^{ i (k+1) \alpha}\: ; \:j \ge 1\Big\}. 
$$
We now derive from \eqref{HS-1500*} that 
\begin{multline*}
\mbox{trace} \Big(M_t^{-1}\bT_{\alpha,t} M_t (M_t^{-1}\bT_{\alpha,t}M_t )_{-2t^{k+1} e^{i (k+1) \alpha}} \Big) \\[6pt]
= \sum_{j} \frac{1}{(\mu_j^{k+1} - t^{k+1} e^{ i (k+1) \alpha})(\mu_j^{k+1} + t^{k+1} e^{ i (k+1) \alpha})},  
\end{multline*}
which yields, since $\alpha = \frac{\pi}{4 (k+1)}$, 
\begin{equation}
\label{HS-1500**}
\mbox{trace}\Big(M_t^{-1}\bT_{\alpha,t} M_t (M_t^{-1}\bT_{\alpha,t}M_t )_{-2t^{k+1} e^{i (k+1) \alpha}} \Big)  = \sum_{j} \frac{1}{\mu_j^{2(k+1)}  - i t^{2(k+1)}}. 
\end{equation}

We have,  by \Cref{prop-pre}, 
\[
\limsup_{|\mu_j| \to + \infty}\left | \frac{\Im (\mu_j)}{\mu_j} \right |  =0.
\]
As a consequence and as in \cite[Proof of Corollary 3]{MinhHung2}, we derive that
\begin{equation}
\label{HS-151}
\sum_j  \frac{1}{\mu_j^{2k+2} - i t^{2k+2}} - \sum_j \frac{1}{|\mu_j|^{2k+2} - i t^{2k+2}} = o(t^{ - 2k - 2+\frac{d}{2}}) \mbox{ as }t \to + \infty.   
\end{equation}
Combining \eqref{HS-1500**} and \eqref{HS-151} yield 
\begin{multline}
\label{HS-p1}
\mbox{trace}\Big(M_t^{-1}\bT_{\alpha,t} M_t (M_t^{-1}\bT_{\alpha,t}M_t )_{-2t^{k+1} e^{i (k+1) \alpha}} \Big)  \\[6pt]
= \sum_{j} \frac{1}{|\mu_j|^{2(k+1)}  - i t^{2(k+1)}}  + o(t^{ - 2k - 2+\frac{d}{2}}) \mbox{ as }t \to + \infty. 
\end{multline}

Applying \eqref{HS-Rm} with $\gamma_1 = t^{k+1} e^{i (k+1) \alpha}$ and $\gamma_2 = -2t^{k+1} e^{i (k+1) \alpha}$ and  using \Cref{lem-mod-res}, we derive that 
\begin{equation}
\label{HS-p111}
(M_t^{-1}\bT_{\alpha,t}M_t )_{-2t^{k+1} e^{i (k+1)  \alpha}} = M_t^{-1}\bT_{\beta, t} M_t. 
\end{equation}
Since 
\[
\mbox{trace}\left (M_t^{-1} \bT_{\alpha, t}\bT_{\beta, t} M_t\right ) = \mbox{trace} \left ( \bT_{\alpha, t}\bT_{\beta, t}\right ), 
\]
it follows from \eqref{HS-p1} and \eqref{HS-p111} that 
\begin{equation}
\label{HS-p2}
\mbox{trace}  \left ( \bT_{\alpha, t}\bT_{\beta, t}\right )
= \sum_{j} \frac{1}{|\mu_j|^{2(k+1)}  - i t^{2(k+1)}} +  o(t^{ - 2k - 2+\frac{d}{2}}) \mbox{ as }t \to + \infty. 
\end{equation}

Applying \Cref{prop-ker}, we derive from \eqref{HS-p2} that, as $t \to + \infty$
\begin{equation}
\label{HS-271}
\sum_j \frac{1}{|\mu_j|^{2k+2}- i t^{2k+2}}  = {\bf c}t^{-2k-2+\frac{d}{2}}+ o(t^{-2k-2+\frac{d}{2}}).
\end{equation}
Considering the imaginary part of \eqref{HS-271} we get, for $\tau = t^{4k+4}$,
\[
\sum_j \frac{1}{|\mu_j|^{4k+4}+\tau} = \Im({\bf c}) \tau^{\frac{d}{8k+8}-1}+ o(\tau^{\frac{d}{8k+8}-1})\mbox{ as }\tau \to + \infty.
\]
Since $\lambda_j = \mu_j+\lambda^*$, it follows that, as $\tau \to + \infty$,
\begin{equation}\label{HS-p4}
\sum_j \frac{1}{|\lambda_j|^{4k+4}+\tau} = \Im({\bf c}) \tau^{\frac{d}{8k+8}-1}+ o(\tau^{\frac{d}{8k+8}-1}).
\end{equation}
Using the fact 
\begin{equation*}
\sum_j \frac{1}{|\lambda_j|^{4k+4}+\tau} = \int_0^{\infty} \frac{d\mathcal{N}(s^{\frac{1}{4(k+1)}})}{s+\tau}, 
\end{equation*}
we derive that 
\begin{equation}\label{HS-p5}
\int_0^{\infty} \frac{d\mathcal{N}(s^{\frac{1}{4(k+1)}})}{s+\tau}  = \Im({\bf c}) \tau^{\frac{d}{8k+8}-1}+ o(\tau^{\frac{d}{8k+8}-1}) \mbox{ as }\tau \to + \infty.
\end{equation}
Applying a Tauberian Theorem of Hardy and Littlewood (see,  e.g., \cite[Theorem 14.5]{Agmon}), we obtain
\[
\mathcal{N}(t) = \frac{\Im ({\bf c})}{\frac{d}{8(k+1)} \int_0^\infty s^{\frac{d}{8(k+1)}-1}(1+s)^{-1}ds} t^{\frac{d}{2}} + o(t^{\frac{d}{2}})\quad \quad  \mbox{ as }t \to + \infty, 
\]
which is the conclusion.

It remains to prove  \eqref{prop-spec-ct}. 
Applying \eqref{HS-Rm} with $\gamma_1 = t^{k+1} e^{i (k+1) \alpha}$ and $\gamma_2 = -2(t+s)^{k+1} e^{i (k+1) \alpha}$ and  using \Cref{lem-mod-res}, we derive that 
$$
(M_t^{-1}\bT_{\alpha,t} M_t )_{-2(t+s)^{k+1} e^{i (k+1)  \alpha}} = M_r^{-1}\bT_{\widetilde{\alpha}, r} M_r. 
$$
where 
$$
\widetilde{\alpha} = \alpha +\frac{\pi}{k+1} \quad \mbox{ and } \quad r = (2(t+s)^{k+1} - t^{k+1})^{\frac{1}{k+1}}. 
$$
Thus by \cite[Theorem 12.14]{Agmon}, 
\begin{equation}
\label{chara-op-2}
\Big \{
s_j +2  (t+s)^{k+1}e^{i\alpha (k+1)}~;~ j\geq 1 
\Big \} \mbox{ is the set of characteristic values of }M_r^{-1}\bT_{\widetilde{\alpha}, r} M_r.
\end{equation}

Applying \Cref{cor-p} and using \eqref{def-Mt}, we have 
\begin{equation}
\label{ct-p-0}
\vvvert M_{t}^{-1} \bT_{\alpha, t} M_{t} \vvvert \leq C t^{-k + \frac{d}{4}} \quad \mbox{ and } \quad \vvvert M_{r}^{-1} \bT_{\widetilde{\alpha}, r} M_{r} \vvvert \leq C r^{-k+\frac{d}{4}}
\end{equation}
for some constant $C>0$ which does not depend on $s$ (and $t$). By \cite[Theorem 12.12]{Agmon} we have
\begin{equation}
\label{ct-p-01}
|\mbox{trace}\Big(M_t^{-1}\bT_{\alpha,t} M_t (M_t^{-1}\bT_{\alpha,t}M_t )_{-2\alpha^{k+1} (t+s)^{k+1}} \Big)| \leq   \vvvert M_{t}^{-1} \bT_{\alpha, t} M_{t} \vvvert  \vvvert M_{r}^{-1} \bT_{\widetilde{\alpha} , r} M_{r} \vvvert.
\end{equation}
Since $-k+\frac{d}{4}<0$  it follows from  \eqref{ct-p-0} and \eqref{ct-p-01} that
\begin{equation}\label{HS-p10}
\lim_{s \to + \infty} \mbox{trace}\Big(M_t^{-1}\bT_{\alpha,t} M_t (M_t^{-1}\bT_{\alpha,t}M_t )_{-2\alpha^{k+1} (t+s)^{k+1}} \Big) = 0. 
\end{equation}
On the other hand, 
by \cite[Theorem 12.14]{Agmon}, 
\begin{multline}\label{HS-p11}
\Big|\sum_{j} \frac{1}{s_{j}(s_{j}+2e^{i\alpha (k+1)} (t+s)^{k+1})} \Big|^2  \le \sum_{j} \frac{1}{|s_{j}|^2} \sum_{j} \frac{1}{|s_{j}+2e^{i\alpha (k+1)} (t+s)^{k+1}|^2} \\[6pt]
\mathop{\le}^{\eqref{chara-op-2}} \vvvert M_{t}^{-1} \bT_{\alpha, t} M_{t} \vvvert^2  \vvvert M_{r}^{-1} \bT_{\widetilde{\alpha} , r} M_{r} \vvvert^2 \mathop{\to}^{\eqref{ct-p-0}} 0 \mbox{ as } s \to + \infty.
\end{multline}
Combining \eqref{HS-p10} and \eqref{HS-p11} yields $c_{t} = 0$,  which is \eqref{prop-spec-ct} .

\medskip 
The proof is complete. 
\end{proof}

\subsection{Proof of \Cref{thm1}} As in \cite[p.34]{MinhHung2}, we derive from  \Cref{prop-ker} that 
\[
\begin{aligned}
 \Im ({\bf c}) &= \frac{1}{(2\pi)^d} \int_{\Omega} \sum_{\ell=1}^2 \int_{\R^d}\frac{1}{\left (  \Sigma_{\ell} (x_0)^{-1} A(x_0) \xi \cdot \xi \right )^{4k+4} +1 } d \xi dx \\
&= \frac{1}{(2\pi)^d}  \sum_{\ell=1}^2\int_\Omega\left | \left \{ \xi : A(x) \xi \cdot \xi < \Sigma_{\ell} (x)\right \}\right | dx \frac{d}{8(k+1)}\int_0^\infty s^{\frac{d}{8(k+1)}-1}(1+s)^{-1} ds.
\end{aligned}
\]
The conclusion now follows from \Cref{prop-spec}. \qed

\section{Completeness of the generalized eigenfunctions of the transmission eigenvalue problem - Proof of \Cref{thm2}} \label{sect-C} 

By \Cref{lem-mod-res}, for all $\widetilde{\theta} \in (0, 2 \pi) \setminus \{\pi \}$, there exists  $t_{\widetilde{\theta}}>0$ such that,  for  $t>t_{\widetilde{\theta}}$,   
\begin{equation}
\label{compl-20}
(T_{\lambda^*}^{k+1})_{t e^{i\widetilde{\theta}} }  = M_{t_k}^{-1}\bT_{\theta,t_k} M_{t_k}, 
\end{equation}
where
\[
\theta = \frac{\widetilde{\theta}}{k+1}  \quad \mbox{ and }\quad t_k = t^{\frac{1}{k+1}}.
\]
By \Cref{prop-pre} and \Cref{cor-p},
\begin{equation}
\label{compl-3-1}
\vvvert  M_{t_k}^{-1}\bT_{\theta,t_k} M_{t_k}  \vvvert \leq Ct_k^{-k+\frac{d}{4}}\quad \quad  \mbox{ and } \quad \quad  \|  M_{t_k}^{-1}\bT_{\theta,t_k} M_{t_k} \|_{L^2(\Omega) \to L^2(\Omega)}\leq Ct_k^{-k}. 
\end{equation}
In particular, $T_{\lambda^*}^{k+1}$ is a Hilbert-Schmidt operator;  moreover, for $t > t_\theta$,   
\begin{equation}
\label{compl-3-2}
\vvvert (T_{\lambda^*}^{k+1})_{t e^{i\widetilde{\theta}}} \vvvert \leq C_{\widetilde{\theta}} t^{-1+\frac{1}{k+1}+\frac{d}{4(k+1)}}.
\end{equation}
Since $k = [d/2] + 1$, it follows that $-1+\frac{1}{k+1}+\frac{d}{4(k+1)} \leq 0$. This implies  that 
\begin{equation}
\label{compl-3-20}
\mbox{for all }\widetilde{\theta} \in (0, 2 \pi) \setminus \{\pi \} \mbox{ there exist }t_{\widetilde{\theta}}>0 \mbox{ and }C_{\widetilde{\theta}}>0 \mbox{ such that }\sup_{t>t_{\widetilde{\theta}}} \vvvert   (T_{\lambda^*}^{k+1})_{t e^{i\widetilde{\theta}}} \vvvert \leq C_{\widetilde{\theta}}.
\end{equation}

Since $T_{\lambda^*}^{k+1}$ is a Hilbert-Schmidt operator, it follows from \cite[Theorem 16.4]{Agmon} that
\begin{itemize}
\item[i)] the space spanned by the generalized eigenfunctions of $T_{\lambda^*}^{k+1}$ is equal to $\overline{\mbox{range}(T_{\lambda^*}^{k+1})}$, the closure of the range of $T_{\lambda^*}^{k+1}$ with respect to the $L^2$-topology.
\end{itemize}
In fact, in order to be able to apply  \cite[Theorem 16.4]{Agmon}, one requires the assumptions on the directions of the minimal growth of the modified resolvent of $T_{\lambda^*}^{k+1}$. We have only proved \eqref{compl-3-2} and \eqref{compl-3-20} instead of this requirement.  Nevertheless, this is sufficient to derive 1) using almost the same arguments in \cite{Agmon} (see also \cite{Robbiano16}). 

\medskip 
The rest of the proof is as in \cite{MinhHung2, Fornerod-Ng1}.  We have
\begin{itemize}

\item[ii)] $\mbox{range } T_{\lambda^*}^{k+1}$ is dense in $[L^2(\Omega)]^2$ since $\mbox{range }T_{\lambda^*}$ is dense in $[L^2(\Omega)]^2$ and $T_{\lambda^*}$ is continuous,

\item[iii)] the space spanned by the general eigenfunctions of $T_{\lambda^*}^{k+1}$ associated to the non-zero eigenvalues of $T_{\lambda^*}^{k+1}$ is equal to the space spanned by the general eigenfunctions of
$T_{\lambda^*}$ associated to the non-zero eigenvalues of $T_{\lambda^*}$ . This can be done as in the last part of the proof of \cite[Theorem 16.5]{Agmon}. Consequently, the space spanned by all generalized
eigenfunctions of $T_{\lambda^*}^{k+1}$ is equal to the space spanned by all generalized eigenfunctions of
$T_{\lambda^*}$.

\end{itemize}
The conclusion now follows from i), ii), and iii). \qed

\bigskip
\noindent {\bf Acknowledgement .} The authors thank Fioralba Cakoni for attracting their attention to the problem and stimulating discussions.  

\providecommand{\bysame}{\leavevmode\hbox to3em{\hrulefill}\thinspace}
\providecommand{\MR}{\relax\ifhmode\unskip\space\fi MR }
\providecommand{\MRhref}[2]{%
  \href{http://www.ams.org/mathscinet-getitem?mr=#1}{#2}
}
\providecommand{\href}[2]{#2}

\end{document}